\documentclass[a4 paper, 12pt]{article}

\usepackage{amsfonts}
\usepackage{amssymb}
\usepackage{amsmath}
\usepackage{fancyhdr}
\usepackage{eufrak}
\usepackage{latexsym}
\usepackage{amsbsy}
\usepackage[all]{xy}
\usepackage[latin1]{inputenc}
\usepackage[T1]{fontenc}        
\usepackage{amsmath}

\newenvironment{proposition}[1]
{\vskip 2mm\noindent \textbf{Proposition #1}  \it}{\vskip 2mm}

\newtheorem{defn}{Definition}[section]
\newtheorem{thm}[defn]{Theorem}
\newtheorem{prop}[defn]{Proposition}
\newtheorem{lem}[defn]{Lemma}
\newtheorem{cor}[defn]{Corollary}
\newtheorem{rem}[defn]{Remark}

\newcommand{\proof}{\vskip 2mm \noindent {\textsc{Proof : }}\rm}
\newcommand{\proo}{\vskip 2mm \noindent {\textsc{Proof }}\rm}
\newcommand{\fin}{\hfill{\Large$\Box$}\\}

\newcommand{\al}{\alpha}

\newcommand{\ga}{\gamma}

\newcommand{\si}{\sigma}

\newcommand{\om}{\omega}\newcommand{\epsi}{\epsilon}

\newcommand{\C}{\mathbb {C}}

\newcommand{\N}{\mathbb {N}}
\newcommand{\Z}{\mathbb {Z}}
\newcommand{\B}{\mathbb {B}}
\newcommand{\U}{\mathbb {U}}
\newcommand{\Pj}{\mathbb {P}}

\newcommand{\orb}{\textsf {O}}
\newcommand{\cd}{\textsf {Card}}
\newcommand{\Id}{{\rm Id}}

\newcommand{\Lip}{{\rm Lip \, }}
\newcommand{\Jac}{{\rm Jac \,  }}

\def\abs#1{\vert #1\vert}
\def\norm#1{\left\|\, #1\,\right\|}

\def\RRR{{\mathfrak R}}
\def\TT{{\cal T}}
\def\Rep{{{\mathbf{R}}}}
\def\Fix{{{\mathbf{F}}}}
\def\EE{{\cal E}}

\def\AA{{\cal A}}
\def\XX{{\cal X}}

\def\CC{{\cal C}}

\def\DD{{\cal D}}
\def\KK{{\cal K}}
\def\FF{{\cal F}}
\def\UU{{\cal U}}
\def\GG{{\cal G}}

\def\RR{{\cal R}}
\def\II{{\cal I}}
\def\LL{{\cal L}}
\def\VV{{\cal V}}
\def\MM{{\cal M}}
\def\GG{{\cal G}}
\def\SS{{\cal S}}
\def\OO{{\cal O}}
\def\QQ{{\cal Q}}
\def\NN{{\cal N}}
\def\WW{{\cal W}}

\def\HH{{\cal H}}
\def\com{\ar@{}[rd]|{\circlearrowleft}}

\title {Normalization of bundle holomorphic contractions and applications to dynamics}
\author{F. Berteloot, C. Dupont, L. Molino}
\date{ \today }

\begin{document}

\maketitle
 
\begin{abstract}
We establish a Poincar\'e-Dulac theorem for sequences $(G_n)_{n \in \Z}$ of holomorphic 
contractions whose differentials $d_0 G_n$ split regularly.
 The resonant relations determining the normal forms hold on the moduli of the exponential rates of contraction. Our results are actually stated in the framework of bundle maps.
 
 Such sequences of holomorphic contractions appear naturally as iterated 
 inverse branches of endomorphisms of $\C \Pj^k$. In this context, our normalization
 result allows to  estimate precisely the distortions of ellipsoids along
 typical orbits. As an application, we show how the Lyapunov exponents of the 
 equilibrium measure are approximated in terms of the multipliers of the repulsive 
 cycles.
\end{abstract}

\section{Introduction and results}\label{intro}

As it is well known, any holomorphic function $F$ which is invertible and 
contracting at the origin of $\C$ is conjugated to its linear part
$A:=F'(0)$. Moreover, the conjugation 
is realized by a function $N$ obtained by a renormalization procedure
: $N=\lim_n A^{-n}F^n$. Remarkably, the proof only relies on the fact that,
by analyticity, $F$ is tangent to $A$ at the order 2.
It turns out that these arguments may be adapted to the case
of holomorphic mappings and yield to the Poincar\'e-Dulac theorem.
We show in this paper that this strategy also works in a non-autonomous
setting, that is for families of holomorphic contractions.\\
 
Our results will be 
expressed in the framework of bundle maps which we now briefly describe (see
the subsection \ref{tube} for more details).
Let $X$ be a set and $\XX :=
X \times \C^k$. For any function $r$ on $X$, the \emph{tube} $E(r) \subset \XX$ is defined by 
\[ E(r) := \bigcup_{x \in X} \{(x,v) \in \XX \ , \ \abs{v} < r(x) \},\]  
where $\abs{\, . \,}$ is the euclidean norm on $\C^k$. 
A \emph{bundle map} over a bijective map $\tau : X \to X$ is a map of the form

\begin{displaymath} 
 \KK :
\begin{array}{rccc}
 &  E(r)          &  \longrightarrow    &         \XX   \\
 &  ( x , v)   & \longmapsto &  \left( \tau( x), K_{ x}(v) \right)
\end{array}                                      
\end{displaymath}

\noindent where $K_x$ is holomorphic and $K_x(0)=0$. A bundle map $\KK$ is
\emph{tame} when the coefficients of the Taylor expansion of $K_{\tau^n(x)}$ grow exponentially slowly 
when $n$ tends to infinity. A function $\phi_\epsi : X \to ]0,1]$ is
    \emph{$\epsi$-slow} if $\phi_\epsi(\tau^{\pm 1}(x)) \geq e^{-
      \epsilon} \phi_\epsi (x)$. We note
$\II$ the identity map of $\XX$ over $\Id_X$. \\

Let us now present our results. To this purpose we will describe the autonomous
case (Poincar\'e-Dulac theorem) and, step by step, precise the corresponding 
statements for bundle maps.  
Let $F: \C^k \to \C^k$ be an holomorphic map whose linear part $A$ satisfies 
$m\vert v\vert \le \vert A(v)\vert \le M\vert v\vert$. If $M^{q+1}< m$,
the sequence $(A^{-n}F^n)_n$ converges as soon as $F$ and $A$ are tangent at 
the order $q+1$. This means that
a sufficently high tangency between $F$ and $A$ allows to overcome the fact that
$\Vert A\Vert \ne \Vert A^{-1}\Vert^{-1}$ and implies the conjugation (when $k=1$, then $m=M$ and one sees again that the tangency at 
the order 2 suffices).
Moreover, the linear map $A$ may actually be replaced by any 
automorphism $N$ which is tangent to $F$ at the order $q+1$. The same phenomenon actually occurs in a non-autonomous setting,
the precise statement is the following :

\begin{thm}\label{HOC}
Let $\FF$ and $\NN$ be two tame bundle maps over $\tau$ which are tangent at the
order $q+1$, with $q \geq 1$. Assume that their linear part
$\AA$ satisfies $m \abs{v} \leq \abs{A_x(v)} \leq M  \abs{v}$
for any $x\in X$, where $0 < m \leq M <1$ and $M^{q+1} < m$. Then $\FF$ is conjugated to $\NN$. \\
More precisely, there exists a $\epsi$-slow function $\rho_\epsi$ such
that $\FF$, $\NN$ contract the tubes $E(\rho_\epsi)$,
$E(2\rho_\epsi)$ and 
a bundle map $\TT:=\lim_{n \to +\infty}
\NN^{-n} \FF^n$ which is $\kappa$-tangent to $\II$ at the order
$q+1$ such
that the following diagram commutes :
\[
\xymatrix{
          E(\rho_\epsi) \ar[r]^\FF \ar[d]_\TT & E(\rho_\epsi) \ar[d]^\TT \\
           E(2\rho_\epsi)  \ar[r]^\NN          &  E(2\rho_\epsi)\;.         }
\] 
 Moreover, for any $M^{q+1}/m < \theta < 1$, 
	     the
	    following estimate 
\[ \forall n \geq 1 \ , \ \forall v \in E_x(\rho_\epsi) \ , \ \abs{(N_x^n \circ T_x - F_x^n)(v)} \leq \varphi_\epsi(x) (m\theta)^n \abs{v}^{q+1}\]
 occurs for some $\epsi$-fast function $\varphi_\epsi$. 
\end{thm}

Let us return to the autonomous case.
To obtain a conjugation between $F$ and $A$ it thus would suffice to perform 
local changes of coordinates cancelling the non-linear terms of order at most 
$q$ in the Taylor expansion of $F$. The determination of these changes of
coordinates is a purely algebraic problem which yields to the so-called 
Poincar\'e's homological equations. As one easily sees, these equations 
have solutions when there are no resonance relations among the eigenvalues of 
$A$. In that case $F$ is thus linearizable. When resonances do occur, 
a finite number of monomials in the Taylor expansion of $F$ can not be cancelled.
However, as it turns out, these terms added to $A$ define a triangular automorphism 
$N$. Then, the sequence $(N^{-n}F^n)_n$ converges and the map $F$ is conjugated to $N$. This is the
Poincar\'e-Dulac theorem for an holomorphic contraction.\\  

This procedure may be used for a bundle map whose linear part is regular and contracting, a property which we state here formally :

\begin{defn}\label{rc}
A linear bundle map $\AA$ over $\tau$ is regular contracting if there exist 

- an integer $1 \leq l \leq k$ and a decomposition $k = k_1 +
\ldots + k_l$,
 
- real numbers $\Lambda_l < \ldots < \Lambda_1 <0$ and $\epsi \ll \abs{\Lambda_1}$ such that :
\[ \ \AA(\LL^j) = \LL^j \ \ \textrm{and} \ \  e^{\Lambda_j - \epsilon} \abs{v}
\leq \abs{A_x(v)}  \leq e^{\Lambda_j + \epsilon}  \abs{v} \;\textrm{for any}\; (x,v) \in \LL^j \;\textrm{and}\; 1\leq j \leq l, \]
where $\LL^j :=  X \times \left[ \{ 0 \} \times \ldots \times \C^{k_j}  \times \ldots \times \{0 \} \right]$ so that $\XX = \LL^1 \oplus \ldots \oplus \LL^l$.

\end{defn}

We say that a bundle map is regular contracting if its linear part is.
In this setting, the resonances hold on the moduli of the contraction rates of $\AA$, that is the $\Lambda_j$.
These are relations of the form $\alpha\cdot\lambda=\Lambda_j$ where $\alpha:=(\alpha_1,\ldots,\alpha_k)\in \N^k$
and $\lambda:=(\Lambda_1 , \ldots , \Lambda_1, \cdots ,  \Lambda_j , \ldots ,
\Lambda_j ,\cdots,
\Lambda_l, \ldots, \Lambda_l)$ (see the subsection \ref{bresonn} for more details). As in the autonomous case, the monomials
$z_1^{\alpha_1}\cdot\cdot\cdot z_k^{\alpha_k}$ in the $j$-th component of the maps $F_x$ may be cancelled by suitable 
changes of coordinates if $\alpha\cdot\lambda\ne\Lambda_j$. This is the meaning of the proposition \ref{poincare} which,
combined with the above theorem \ref{HOC} gives the following version of Poincar\'e-Dulac theorem for bundle maps :

\begin{thm}\label{normaliza}
A tame and regular contracting bundle map $\GG$ over $\tau$ is conjugated to a resonant bundle map $\RR$. More precisely there exist a bundle map $\VV$ 
which is $\kappa$-tangent to $\II$ and a $\epsi$-slow function $r_\epsi$ such that the following diagram commutes :
\[
\xymatrix{
          E(r_\epsi) \ar[r]^\GG \ar[d]_\VV   & E(r_\epsi) \ar[d]^\VV \\
           E(2r_\epsi)  \ar[r]^{\RR}         &  E(2r_\epsi)        }
\]
and $\GG,\RR$ contract the tubes $E(r_\epsi), E(2r_\epsi)$.
\end{thm}

The above theorem is a linearization statement ($\RR=\AA$) when there are no resonances among the $\Lambda_j$. Let us however stress that even in the resonant case,
the ``stability''
properties of resonant bundle maps (see propositions \ref{olk} and \ref{pj}) imply that the iterated bundle maps $\RR^n$ and $\AA^n$ behave similarly, a fact which is of
great importance for our applications.\\

Similar results were proved by Guysinsky-Katok-Spatzier (\cite{KS}, \cite{GK} theorem 1.2) for smooth bundle maps and by Jonsson-Varolin (\cite{JV} theorem 2) 
in the holomorphic case. In this article, Jonsson and Varolin proved the
Bedford's conjecture on the complex structure of stable manifolds of
holomorphic automorphisms in the non-uniform setting. The articles of
Forn\ae ss-Stens\o nes \cite{FS} and Peters \cite{P} are also dedicated to this conjecture.\\

The originality of our 
approach is the use of a renormalization technique. It has the advantage of being simpler and 
also provides an answer to a question asked by Jonsson and Varolin (\cite{JV} final remarks). To our knowledge, this approach also gives the simplest proof of Poincar\'e-Dulac theorem for an holomorphic contraction (see the survey \cite{B} for a precise exposition). Other simple proofs are due to Sternberg \cite{St} and Rosay-Rudin \cite{RR}.\\

The Oseledec-Pesin reduction theorem (see theorem \ref{OPR}) opens a large field of applications for the theorem \ref{normaliza} in the setting of smooth 
ergodic dynamical systems. We will here investigate the case of holomorphic endomorphisms of the complex projective space $\C\Pj^k$. Our aim is to precisely describe the asymptotic behaviour of typical iterated inverse branches for such endomorphisms.
Let us recall that the works of Fornaess-Sibony \cite{S} and Briend-Duval \cite{BD} show that any holomorphic endomorphism $f:\Pj^k\to\Pj^k$ of algebraic degree 
$d \ge 2$ induces an ergodic dynamical system $(\Pj^k,f,\mu)$ where $\mu$ is the unique maximal entropy measure of $f$. This measure is mixing and its Lyapunov exponents  
$\chi_1 \leq \ldots \leq \chi_k$ are bounded from below
by $\log \sqrt d$.\\
We aim to apply the Oseledec-Pesin reduction theorem and the normalization theorem \ref{normaliza} to a bundle map generated by the inverse branches of $f$.
For that purpose we work in the set of orbits
$\orb := \{ \hat x := (x_n)_{n \in \Z} \, , \,  x_{n+1} = f (x_n)
\}$. The right shift $\tau$ acts on $\orb$ and leaves invariant a probability measure 
 $\nu$ related to $\mu$ by $\nu (\pi ^{-1}
(A)) = \mu(A)$, $\pi$ being the time zero projection $\pi(\hat x)=x_0$.
 A typical orbit $\hat x$ does not intersect the critical locus of $f$ for all $n$ and we may therefore define    
 the inverse branch $f^{-n}_{\hat x}$ of
  $f^n$ that sends $x_0$ to $x_{-n}$. Our result compares $f^{-n}_{\hat x}$ with its linear tangent map : 

\begin{thm}\label{nt}
Let $f$ be an holomorphic endomorphism of algebraic degree $d$ on $\Pj^k$. Let $\Sigma_s$ be the sum of the $s$ largest Lyapounov exponents of the
maximal entropy measure $\mu$. Let $\nu$ be the measure induced by $\mu$ on the natural extension $\orb$ of $(\Pj^k,f,\mu)$ and $\tau$ the right
 shift on $\orb$. 
There exist a full measure subset $X$ of $\orb$, $\epsi$-slow functions $r_\epsi, t_\epsi : X \to ]0,1]$, a
resonant bundle map $\RR$ over $\tau$ and an injective bundle map $\SS$ over $\Id_X$ such that the following diagram commutes for 
all $n \geq 1$ (we note $R^n_{\hat x } = R_{\tau^n(\hat x)} \circ \ldots \circ  R_{\hat x}$) :
\[  
\xymatrix{
    B_{x_0}(r_\epsi(\hat x)) \ar[rr]^{f^{-n}_{\hat x}} \ar[d]_{S_{\hat x}}  & &
      f^{-n}_{\hat x} [ B_{x_0}(r_\epsi(\hat x)) ] \ar[d]^{S_{\tau^n(\hat x)}}   \\
      \B(t_\epsi (\hat x)) \ar[rr]^{R^n_{\hat x}}    &   &   \B(t_\epsi (\tau^n(\hat x)))  }
  \]
There exist also constants $\al,M >0$ and $\epsi$-fast functions $\beta_\epsi, L_\epsi,T_\epsi : X \to [1,+\infty[$ such that for all $n \geq 0$ : 
\begin{enumerate}
\item $f^{-n}_{\hat x} [ B_{x_0}(r_\epsi(\hat x)) ] \subset B_{x_{-n}}(M)$,
\item $\forall (p,q) \in   f^{-n}_{\hat x} [B_{x_0}(r_\epsi(\hat x)) ] \, , \, \al \, d(p,q) \leq 
\abs{S_{\tau^n(\hat x) }(p) - S_{\tau^n(\hat x)}(q) } \leq  \beta_{\epsi}(\tau^n(\hat x)) \, d(p,q)$,
\item $\Lip f_{\hat x}^{-n} \leq L_\epsi(\hat x) e^{-n \chi_1 + n \epsi}$ on $B_{x_0}(r_\epsi(\hat x))$,
\item for all $p \in f_{\hat x}^{-n}[B_{x_0}(r_\epsi(\hat x))] , \, \vert {1 \over n} \log \norm {\bigwedge^s d_p f^n} - 
\Sigma_s  \vert  \leq {1 \over n } \log T_\epsi(\hat x) + \epsilon$.
\end{enumerate}
\end{thm}

The above result is useful for studying the properties of the maximal entropy measure $\mu$.
A weak version of it, corresponding to the case where all Lyapunov exponents are equal, has been used 
in \cite{BDu} for characterizing the endomorphisms $f$ for which $\mu$ is absolutely continuous with respect to the Lebesgue measure.
In this article, we will use
the theorem \ref{nt} for proving 
the following approximation formula for Lyapunov exponents :

\begin{thm}\label{appli}
We use the same notations than in the previous theorem, recall that $\Sigma_s = \chi_{k-s+1} + \ldots + \chi_k$.
Let $\Rep_n$ (resp. $\Rep_n^*$) be the set of repulsive periodic points whose period divides $n$
(resp. equals $n$). Then : $\Sigma_s =\lim_{n \to
    + \infty}  {1 \over d_t^n} \sum_{p \in \Rep_n} {1 \over n} \log
  \norm{\bigwedge^s d_p f^n} = 
\lim_{n \to
    + \infty}  {1 \over d_t^n} \sum_{p \in \Rep_n^*} {1 \over n} \log
  \norm{\bigwedge^s d_p f^n}$.
\end{thm}

Bedford-Lyubich-Smillie \cite{BLS} 
proved a similar result for the positive Lyapunov exponent of a Hénon map $f : \C^2 \to \C^2$, with $\Rep_n$ replaced 
by the $n$-periodic saddle points. The theorem \ref{appli} was proved by Szpiro-Tucker for rational maps ($k=1$) whose coefficients are in a number field (\cite{ST}, corollary 6.1). Observe also that for $s = k$, the exterior product
$\bigwedge^k d_p f^n$ is the jacobian of $f^n$ at $p$, 
which satisfies the multiplicative property $\bigwedge^k d_p f^{m+n} = \bigwedge^k d_{f^n(p)} f^m . \bigwedge^k d_{p} f^n$. We thus also have :
\begin{cor}
$\lim_{n \to
    + \infty}  {1 \over d_t^n} \sum_{p \in \Rep_n} \log
  \abs{\Jac d_p f} = \chi_1 + \ldots + \chi_k$.
\end{cor}

These approximation formulas have some importance in the study of bifurcations of holomorphic families of endomorphisms of $\Pj^k$.
In particular the theorem 2.2 of Bassanelli-Berteloot \cite{BB} is a consequence of the above corollary.\\

We may extend the theorems \ref{nt} and \ref{appli} to polynomial-like mappings i.e. holomorphic and proper maps $F : U \to V$ between two open sets 
$U \Subset V \subset \C^k$.  The dynamical properties of these maps have been studied by Dinh-Sibony \cite{DS}. 
They proved in particular the existence of an equilibrium measure $\mu$ when $d_t\ge 2$. This measure is mixing and its Lyapunov exponents are positive if $\mu$ is PLB (i.e. the plurisubharmonic functions are in 
$L^1(\mu)$). The theorem \ref{appli} remains valid if the cardinal of the $n$-periodic points is asymptotically bounded from above by $d_t^n$ (\cite{DS} subsection 3.5). 
These conditions are fullfilled e.g. for perturbations of polynomial lifts of endomorphisms of $\Pj^k$.

\section{Generalities}

\subsection{Slow and fast functions, tube and bundle maps}\label{tube}

We introduce in this subsection several notations and definitions. We consider a set $X$ and a fixed
bijective map $\tau : X \to X$. A function $\phi_\epsi : X \to ]0,1]$ is
    \emph{$\epsi$-slow} if $\phi_\epsi(\tau^{\pm 1}(x)) \geq e^{-
      \epsilon} \phi_\epsi (x)$. A function $\phi : X \to ]0,1]$ is
    \emph{slow} if $\phi$ is bounded from below by a $\epsi$-slow
    function for all $\epsi \ll 1$. A function $\psi_\epsi : X \to [1,+\infty[$ is \emph{$\epsi$-fast} if $1/\psi_\epsi$ is $\epsi$-slow and a function $\psi : X \to [1,+\infty[$ is fast if $1/\psi$ is slow.

\begin{lem}\label{abelard} Assume that $X$ is endowed with a $\tau$-invariant probability measure
 $\nu$. Let $\epsi > 0$ and $u : X \to ]0 , 1]$ satisfying $\log u \in L^1(\nu)$. Then there
	    exists a $\epsi$-slow function $u_\epsi$ such that $u_\epsi \leq u$ for $\nu$-almost every $x \in X$. 

Similarly, if $v : X \to [1 ,
	  +\infty [$ satisfies $ \log v \in L^1(\nu)$ then there exists a
	    $\epsi$-fast function $v_\epsi$ such that $v \leq v_\epsi$ for $\nu$-almost every $x \in X$. 
\end{lem}

\proof By Birkhoff ergodic theorem, we have $\lim_{n \to \pm \infty} {1 \over n} \log u(\tau^{n}(x)) =  0$ for $\nu$-almost every $x \in X$. So
there exists $V : X \to ]0,1]$ such that $u(\tau^n(x)) \geq e^{-\abs
	  {n} \epsi} V(x)$ for all $n \in \Z$. Let $u_\epsi(x) := \inf_{n \in \Z} \{ u(\tau^n(x)) e^{\abs{n} \epsi}  \}$. We have $V(x) \leq u_\epsi (x) \leq
	u(x)$ and 
\begin{displaymath}
\begin{array}{rcl}
  u_\epsi (\tau (x))  &  =   &  \inf_{n \in \Z}  \{ u(\tau^{n+1}(x)) e^{\abs {n} \epsi}
\}  \\
       &  =   &  e^{-\epsi} \inf_{n \in \Z} \{ \ldots , u(\tau^{-1}(x))
e^{3\epsi} , u(x) e^{2\epsi} ,  u(\tau(x)) e^{\epsi} ,  u(\tau^2(x))
e^{2\epsi} , \ldots   \}   \\
       &  \geq   &  e^{-\epsi} \inf_{n \in \Z} \{ \ldots , u(\tau^{-1}(x))
e^{\epsi} , u(x)  ,  u(\tau(x)) e^{\epsi} ,  u(\tau^2(x))
e^{2\epsi} , \ldots   \}   \\ 
 &  =  & e^{-\epsi} u_\epsi(x). 
\end{array}                                      
\end{displaymath} 
We prove similarly  $u_\epsi (\tau^{-1} (x)) \geq e^{-\epsi}
u_\epsi(x)$ : the function $u_\epsi$ is therefore $\epsi$-slow. The analogous property for $v$ is obtained by considering $u := 1/v$. \fin

We set $\XX = X \times \C^k$. Let $E_x =  \{(x,v) \in \XX \ , \ v \in \C^k  \}$ and $E_x(t) := \{(x,v) \in \XX \ , \ \abs{v} < t  \}$. For any function $r : X \to
  [0,+\infty[$, we note $E_x(r) := E_x(r(x))$ and $E(r) = \bigcup_{x
  \in X} E_x(r(x)) $. The subset $E(r) \subset \XX$ is the \emph{tube}
  of radius $r$. We say that $E(r)$ is \emph{slow} (resp. $\epsi$-slow) if $r$ is slow (resp. $\epsi$-slow). Idem with the ``fast'' terminology. \\

Let $\si \in \{ \Id_X , \tau \}$. A holomorphic bundle map $\KK : E(r) \to
\XX$ over $\si$ is a map satisfying $\KK(x,v) =
(\si(x) , K_x(v))$, where  $K_x : E_x(r) \to E_{\si(x)}$ is
holomorphic and $K_x(0)=0$. We say that $\KK$ is \emph{tame} if there
exists $\epsi_0 > 0$ such that for all $\epsi < \epsi_0$ there exist a $\epsi$-slow function
$r_\epsi$ and a $\epsi$-fast function $s_\epsi$ with  $\KK : E(r_\epsi) \to E(s_\epsi)$.

A \emph{stable} bundle map, that is a map of the form $\KK : E(r) \to E(r)$, may be iterated.
The \emph {n-th iterate} $\KK^n$  is defined by $K^n_x:=K_{\si^{n-1}(x)}\circ\cdot\cdot\cdot\circ K_x$.

 When it makes sense, we consider $\KK^{-1}$ the \emph{inverse} bundle map of $\KK$. We note 
 $\KK^{-1}(x,v) =
(\si^{-1}(x) , K_x^{-1}(v))$ where $K^{-1}_{x} := (
K_{\si^{-1}(x)} )^{-1} : E_x \to E_{\si^{-1}(x)}$ (this map is defined in a neighboorhood of
$0 \in E_x$).\\

We note
$\Lip(\KK) := \sup_{x \in X}  \Lip (K_x)$ and say that $\KK$
\emph{contracts} $E(r)$ if $\Lip(\KK) < 1$ on $E(r)$. The tube $E(r)$
is \emph{stable} if $\KK(E(r)) \subset E(r)$. Let $\II$ be the identity map of $\XX$ over $\Id_X$ and $\OO$ the zero
bundle map over $\tau$, which sends $\XX$ to $X \times \{ 0 \}$. We say that a
bundle map $\KK$ over $\Id_X$ is \emph{$\kappa$-tangent} to $\II$ if
$\Lip(\KK - \II) < \kappa$. \\

Let $m \geq 1$. A bundle map $\KK$ is \emph{homogeneous of degree $m$} (or \emph{$m$-homogeneous}) if  the map $K_x$ is homogeneous of
  degree $m$ for all $x \in X$. Let $\sum_{m
  \geq 1} \KK^{(m)}$ be the Taylor expansion of $\KK$, where
  $\KK^{(m)}$ is $m$-homogeneous. The \emph{linear
  part} of $\KK$ is $\KK^{(1)}$. A bundle map $\KK$ is \emph{polynomial} if there exists $m_0 \geq 1$ such that $\KK = \sum_{m=1}^{m_0} \KK^{(m)}$.  \\

Let $\KK_1$ and $\KK_2$ be two bundle maps over $\si \in \{ \Id_X , \tau \}$. We say that $\KK_1$ is \emph{tangent to $\KK_2$ at the order $m+1$} if $\KK_1$ and
 $\KK_2$ share the same Taylor expansion up to the order
$m$. We then note $\KK_1 = \KK_2 + O(m+1)$. When $m=1$, we just say that $\KK_1$ is
 \emph{tangent} to $\KK_2$. If $\KK_1$ and $\KK_2$ are bundle maps over $\tau$, we say that $\KK_1$ 
is  \emph{conjugated to} $\KK_2$ if there exist a bundle map $\WW$ tangent to $\II$ and positive functions $r,s$ such that the following diagram is commutative :
\[
\xymatrix{
          E(r) \ar[r]^{\KK_1} \ar[d]_\WW & E(r) \ar[d]^\WW \\
           E(s)  \ar[r]^{\KK_2}         &  E(s)         }
\]
We say that $\KK_1$ is \emph{conjugated to} $\KK_2$  \emph{at the order
  $m+1$} if $\KK_1$ is conjugated to a bundle map $\tilde{\KK_2}$ which is tangent to $\KK_2$ at the order $m+1$. \\

Let $\abs{\KK^{(m)}}(x):= \max_{v \neq 0} \abs{K_x^{(m)} (v) } /
  \abs{v}^m$ (also denoted $\abs{K_x^{(m)}}$). For any $\al = (\al_1, \ldots , \al_k) \in \N^k$, we note
  $\abs{\al} := \al_1+ \ldots+\al_k$ and $P_\al (z_1,\ldots,z_k) =
  z_1^{\al_1} \ldots z_k^{\al_k}$. Then $\abs{\KK^{(m)}}(x)$ is equal (up to a constant) to the maximum of
  the coefficients of $K_x^{(m)}$ with respect to the basis $\{ P_\al
  e_i \, , \, \abs {\al} = m \, , \, 1 \leq i \leq k \}$, $(e_i)_{1 \leq i \leq
  k}$ being the canonical basis of $\C^k$. The notation
$\abs{\, . \,}$ will also be used for the standard hermitian norm on the
spaces $(\bigwedge^s \C^k)_{1 \leq s \leq k}$ and for the norm of
  operators : $\bigwedge^s \C^k \to \bigwedge^s \C^k$.

\subsection{Two simple lemmas on bundle maps}

\begin{lem}\label{cow}
Let $\si \in \{ \Id_X, \tau  \}$ and $\KK : E(r)  \to  E(s)$ be a bundle map over $\si$. For all $m \geq 1$, $\abs{\KK^{(m)}}  \leq s \circ \si /r^m$. In particular $\abs{\KK^{(m)}}$ is a fast function when $\KK$ is tame.
\end{lem}

\proof
Let $K_x : E_x(r) \to E_{\si(x)}(s)$. For all $v \neq 0$ and $\rho < r(x)$, we have :
\[  {\rho^m  \over \abs{v}^m } K_x^{(m)}(v) = K_x^{(m)}\left(\rho {v \over \abs{v}}\right) =   {1 \over 2\pi} \int_0^{2\pi} K_x\left(\rho {v \over \abs{v}} e^{i\theta} \right)  e^{-im\theta} \, d\theta . \]
The lemma follows by taking the norm and the limits when $\rho$ tends to $r(x)$.\fin

\begin{lem}\label{NN}
Let $\si \in \{ \Id_X, \tau \}$ and $\KK$ be a tame bundle map over $\si$. Let $\DD$ be the linear
part of $\KK$. Assume that there exist $0< a \leq b$ such that $a
\abs{v} \leq \abs{D_x(v)} \leq  b \abs{v}$. Let $\ga,\kappa > 0$. Then
for $\epsi \ll 1$ there exists a $\epsi$-slow function $\phi_\epsi$ such that : 
\begin{enumerate}
\item $\forall (u,v) \in E_x(\phi_\epsi), a e^{-\ga} \abs{u-v} \leq
  \abs{K_x(u) - K_x(v)} \leq b e^\ga \abs{u-v}$. 

In particular if 
$b e^\ga \leq e^{-\epsi}$ then the tube $E(\phi_\epsi)$ is stable by $\KK$.

\item $\Lip (\KK - \DD) \leq \kappa$ on $E(\phi_\epsi)$.

\item If $\DD = \II$, $\abs{ \bigwedge^s \Id_{\C^k} - \bigwedge^s
  (d_t K_x  )^{\pm 1} } \leq  {1 \over 10}$ for all $t \in E_x(\phi_\epsi)$ and $1 \leq s \leq k$. 
\end{enumerate}
\end{lem}

\proof
Let $\epsi' = \epsi / 3$ with $\epsi \ll 1$. As $\KK$ is tame, there exist a
$\epsi'$-slow function $r_{\epsi'}$ and a $\epsi'$-fast function
$s_{\epsi'}$ such that $\KK : E(r_{\epsi'}) \to E(s_{\epsi'})$. Let $x \in X$. The Cauchy's
estimates on $E_x(r_{\epsi'}/2)$ bound the second derivatives of $K_x$ by $c \, 
s_{\epsi'}(\si(x)) /r_{\epsi'}(x)^2$, where $c$ is a constant depending only on the dimension $k$. We deduce that for all $0 \leq \rho \leq  r_{\epsi'}(x)/2$ :
\begin{equation}\label{aze}
 \forall t \in E_x(\rho) \ , \  \abs{d_t (D_x - K_x)} = \abs{d_0 K_x - d_t K_x} \leq   {c \,  s_{\epsi'}(\si(x)) \over r_{\epsi'}(x)^2} \rho.
\end{equation}

 Let $\eta < 1$ be such that for any $1 \leq s \leq k$ and any linear map $L
: \C^k \to \C^k$, $\abs{ \Id_{\C^k} -L} < \eta$ implies $\abs{  \bigwedge^s \Id_{\C^k} - \bigwedge^s L^{\pm 1} } \leq  {1 \over 10}$. Define the $\epsi$-slow function 
$\phi_\epsi$ by:
\[ \phi_\epsi :=  {r_{\epsi'}^2 \over c \,   s_{\epsi'}\circ \si}  \min \{ (e^\ga - 1) b \, , \,
  (1-e^{-\ga}) a \, , \, \kappa , \eta  \}. \]

\noindent As $\phi_\epsi \leq  {r_{\epsi'}^2 \over c \, 
  s_{\epsi'}\circ \si}$ and $r_{\epsi'} \leq 1 \leq s_{\epsi'}$,
we may assume that $\phi_\epsi \leq r_{\epsi'} /2$ by taking $c \geq
2$. Making $\rho = \phi_\epsi$ in (\ref{aze}) we get the following estimates on $E_x(\phi_\epsi)$ : 
\begin{displaymath}
\begin{array}{rcl}
   \abs {K_x(u) - K_x(v)}   & \leq  & \abs { D_x(u) - D_x(v) } + \abs {(D_x
  - K_x)(u) - (D_x - K_x)(v)}    \\
                   &   \leq   &  b \abs{u-v} +  {c \,  s_{\epsi'}(\si(x)) \over r_{\epsi'}(x)^2}
  \phi_\epsi(x) \abs{u-v} \\
                   & \leq   &  b e^\ga  \abs{u-v}. 
\end{array}                                      
\end{displaymath} 
We have similarly $\abs {K_x(u) - K_x(v)} \geq a e^{-\ga}
\abs{u-v}$. If $b e^\ga \leq e^{-\epsi}$, then $\abs
    {K_x(u)} \leq b e^\ga \phi_\epsi(x) \leq e^{-\epsi} \phi_\epsi(x)
    \leq  \phi_\epsi(\si(x))$ on $E(\phi_\epsi)$, and the tube
    $E(\phi_\epsi)$ is stable by $\KK$. The points $2$ and $3$ are
    also a consequences of (\ref{aze}) with $\rho = \phi_\epsi$, using respectively the estimates $\phi_\epsi \leq  \kappa \,  r_{\epsi'}^2 / c \, s_{\epsi'}\circ \si$  and $\phi_\epsi \leq \eta \, r_{\epsi'}^2 / c \, s_{\epsi'}\circ \si$. \fin

\subsection{Resonances} \label{bresonn}

In this subsection, we will consider a linear bundle map $\AA$ over
$\tau$ which is regular contracting. Our aim is to discuss the resonances associated to such a bundle map.
 According to definition \ref{rc},  
there exist $k = k_1 +
\ldots + k_l$, $\Lambda_l < \ldots < \Lambda_1 <0$ and $\epsi \ll \abs{\Lambda_1}$ such that
\[ \forall 1 \leq j \leq l \ , \  \LL^j :=  X \times \left[ \{ 0 \} \times \ldots \times \C^{k_j}  \times \ldots \times \{0 \} \right],\]
\[ \AA(\LL^j) = \LL^j \ \ \textrm{and} \ \  \forall (x,v) \in \LL^j \ , \ e^{\Lambda_j - \epsilon} \abs{v}
\leq \abs{A_x(v)}  \leq e^{\Lambda_j + \epsilon}  \abs{v}. \]

 \noindent Observe that the matrix of $A_x$ in the canonical basis is block diagonal.  \\

Let us now consider a bundle map $\KK : E(r) \to \XX$ whose linear
part is equal to $\AA$. Let $\pi_j(\KK) : E(r) \to \LL^j$ be the $j$-th component of $\KK$ with 
respect to the splitting $\XX = \oplus_{1 \leq j \leq l} \LL^j $. If $\al \in \N^k$, 
we note $\pi_j^\al (\KK) : \XX \to \LL^j$ the homogeneous part of degree $\al$ in $\pi_j(\KK)$. \\ 

We set $(\lambda_1, \ldots , \lambda_k) :=
(\Lambda_1 , \ldots , \Lambda_1, \cdots ,  \Lambda_j , \ldots ,
\Lambda_j ,\cdots,
\Lambda_l, \ldots, \Lambda_l)$, where $\Lambda_j$ appears $k_j$ times. If $\al = (\al_1,\ldots,\al_k) \in \N^k$, 
we note $\al \cdot \lambda := \al_1\lambda_1 +\ldots +\al_k \lambda_k$. 
We shall denote by  ${\tilde q} \geq 1$ the entire part of $\Lambda_l / \Lambda_1$.

\begin{defn} Let $\AA$ be the above linear bundle map. For any $1 \leq j \leq l$, the set $\RRR_j$ of $j$-\emph{resonant degrees} is defined by : 
\[ \RRR_j := \{ \al \in \N^k \ , \
  \abs{\al} \geq 2 \  \textrm{ and }  \  \al \cdot \lambda =
  \Lambda_j \}. \] 
The set $\mathfrak{B}_j$ of \emph{$j$-subresonant degrees} is
  defined similarly with $\al \cdot \lambda < \Lambda_j$ and the set $\mathfrak{P}_j$ of \emph{$j$-superresonant degrees} with $\al
  \cdot \lambda > \Lambda_j$. We set $\mathfrak{B} := \cup_{j=1}^l \mathfrak{B}_j$, 
 $\RRR := \cup_{j=1}^l \RRR_j$ and $\mathfrak{P} := \cup_{j=1}^l \mathfrak{P}_j$. 
\end{defn}

It should be observed that the sets $\RRR$, $\mathfrak{B}$ and $\mathfrak{P}$ are unchanged if $\AA$ is replaced by $\AA^n$.
As $\Lambda_l < \ldots < \Lambda_1 <0$ one sees that
  $\{ \abs{\al} \geq {\tilde q} + 1 \} \subset \mathfrak{B}$ and $\RRR \cup \mathfrak{P} \subset \{ 2 \leq  \abs{\al} \leq {\tilde q}
  \}$. In particular, $\RRR$ and $\mathfrak{P}$ are finite
  sets. Observe also that the set $\RRR_1$ is empty, and that $\al_i =
  0$ for  $\al \in \RRR_j$ and $i \geq k_1 + \ldots + k_{j-1}+1$. In
  particular for any $j \geq 2$ and $\al \in \RRR_j$, the bundle map $\pi_j^\al(\KK)$ may be
  viewed as a bundle map from $\LL^1 \oplus \ldots \oplus
  \LL^{j-1}$ to $\LL^j$ :
\begin{equation}\label{lok}
  \pi_j^\al(\KK) \, : \, E(r) \cap \left[ \LL^1 \oplus \ldots \oplus
  \LL^{j-1} \right] \longrightarrow \LL^j.  
\end{equation}
The following lemma will be used in subsection \ref{fo}, the proof is left to the reader. We recall that $\epsilon \ll \abs{\Lambda_1}$.
\begin{lem}\label{marge}
There exists $\zeta >0$ such that for all $1 \leq j \leq l$ :
\begin{enumerate}
\item if $\al \in \mathfrak{B}_j$ then $\al \cdot \lambda - \Lambda_j +
  (\abs{\al} +2) \epsilon \leq - \zeta$.
\item if $\al \in \mathfrak{P}_j$ then $\al \cdot \lambda - \Lambda_j - (\abs{\al} +2) \epsilon \geq \zeta$. 
\end{enumerate}
\end{lem}
The disjoint sets $\RRR, \mathfrak{B} , \mathfrak{P}$ lead to the following decomposition for $\KK : E(r) \to \XX$ :
\[  \KK = \KK^{(1)} + \sum_{1 \leq j \leq l \, , \,  \al  \in \RRR_j}  \pi_j^\al (\KK) + \sum_{1 \leq j \leq l \, , \,  \al  \in \mathfrak{B}_j}  \pi_j^\al (\KK) + \sum_{1 \leq j \leq l \, , \,  \al  \in \mathfrak{P}_j}  \pi_j^\al (\KK). \] 
These three sums are respectively denoted $\RRR (\KK)$, $\mathfrak{B}(\KK)$, $\mathfrak{P}(\KK)$ 
and called the \emph{resonant part} and the \emph{sub/superresonant parts} of $\KK$. 
\begin{defn}
A bundle map $\KK$ over $\tau$ is resonant if
$\mathfrak{B}(\KK) = \mathfrak{P}(\KK) = \OO$ (the zero bundle map
over $\tau$). In other words, $\KK$ is resonant if $\KK = \AA + \RRR(\KK)=\KK^{(1)} + \RRR(\KK)$.
\end{defn}

The following classical result will be crucial. It asserts that resonant bundle maps 
enjoy strong stability properties under iteration (see section \ref{pest} for a proof). 

\begin{prop}\label{olk}
Let $\KK$ be a resonant bundle map. For any $n \geq 1$, we have
$\mathfrak{B}(\KK^n) = \mathfrak{P}(\KK^n) = \OO$. In particular $\KK^n =
\AA^n +  \RR(\KK^n)$ and the degree of $\KK^n$ is bounded by ${\tilde q}$. 
\end{prop}

\begin{rem}\label{trian}
Let $\KK$ be a resonant bundle map. It follows from (\ref{lok}) and
the proposition \ref{olk} that for any $x \in X$, $w \in E_x$ and $n \geq 1$, the matrix of the differential $d_w K^n_x$ is ``block lower triangular'', with a block diagonal part equal to $A^n_x$.  
\end{rem}

The following result compares $\abs{\bigwedge^s d_w
    K^n_x}$ and $\abs{\bigwedge^s A^n_x}$ (see section \ref{pest}). 
\begin{prop}\label{estreso}
Let $\KK : E(\rho_\epsi) \to E(\rho_\epsi)$ be a resonant bundle map
    where $\rho_\epsi$ is a $\epsi$-slow function. There exists $\eta > 0$ 
(depending only on $\tilde q, k$ and not on $\KK$) and a $\eta\epsi$-fast function $H_{\eta\epsi} : X \to [1,+\infty[$ such that
    for all $w \in E_x(\rho_\epsi)$, $s \in \{1,\ldots,k \}$ and $n
    \geq 1$ :  
\[  \Big \vert   {1 \over n} \log \abs{\bigwedge^s d_w
    K^n_x} - (\lambda_1 + \ldots + \lambda_s)  \Big \vert  \leq {1
    \over n} \log H_{\eta\epsi}(x) + \eta \epsilon. \]
\end{prop}

\section{Normalization of contracting bundle maps}\label{resonn}

This section is devoted to the proof of theorem \ref{normaliza}.
We essentially proceed in two steps. We first prove that sufficently tangent bundle maps are conjugated. We then show that any bundle map satisfying 
the assumptions of theorem \ref{normaliza} is tangent to a resonant bundle map.

\subsection{High tangency implies conjugation}

The aim of this subsection is to prove the theorem \ref{HOC} stated in the introduction.

\proof
Let $M^{q+1}/m < \theta < 1$ and $\kappa > 0$. We set $\si := \sum_{j
  \geq 0} \theta^j$ and $\beta :=  e^{(q+2)\epsilon} M ^{q+1} / m$ where
$\epsi$ is so small that $M e^\epsi \leq e^{-\epsi}$ and $\beta e^{\epsi} \leq \theta$.

 To establish the theorem, we will show that there exist a $\epsi$-slow
function $r_\epsi$ and a $\epsi$-fast function $D_\epsi$ such that the
map $T_{x,n}:=N^{-n}_{\tau^n(x)} F^n_x$ is well defined on $E_x(r_\epsi)$ and the following estimates hold on $E_x(r_\epsi)$ :
\begin{equation}\label{ineg}
 \max \left\{ \abs{F_x(v)} \ , \  \abs{N_x(v)}
 \right\} \leq M e^\epsi \abs{v},  
\end{equation}
\begin{equation}\label{tyu}
 \abs{  T_{x,n+1}(v) - T_{x,n}(v)  } \leq \beta^n  D_\epsi(\tau^n(x)) \abs{v}^{q+1}.
\end{equation}
Let us note that for any $\epsilon$-slow function $\rho_{\epsilon}\le \frac{r_{\epsilon}}
{2}$ the inequality 
$M e^\epsi \leq e^{-\epsi}$ and (\ref{ineg}) implies that $E(\rho_\epsi)$ and $E(2\rho_\epsi)$ are stable by $\FF$ and $\NN$. \\

We start by showing how the theorem \ref{HOC} follows from (\ref{ineg})
and (\ref{tyu}). We check the convergence of $T_{x,n}$ on $E(r_\epsi)$
and define the functions $\varphi_\epsi$ and $\rho_\epsi \leq
r_\epsi$. The estimate (\ref{tyu}) implies for all $v \in E_x(r_\epsi)$ : 
\begin{equation}\label{sds}
 \abs{  T_{x,n+1}(v) - T_{x,n}(v) } \leq \beta^n D_\epsi(x) e^{n\epsilon}  \abs{v}^{q+1} \leq  \theta^n D_\epsi(x) \abs{v}^{q+1} .
\end{equation}
Therefore $T_x := \lim_{n \to +\infty} T_{x,n}$ exists on
$E_x(r_\epsi)$ and satisfies $T_x =  N_{\tau(x)} ^{-1} T_{\tau(x)}
F_x$. We define $\varphi_\epsi :=  \si D_\epsi$ and $\rho_\epsi := \min \{ \frac{r_\epsi}{2} , \kappa / \varphi_\epsi  \}$. 
Observe that (\ref{sds}) implies on $E_x(\rho_\epsi)$ : 
\[ \abs{ T_x(v) - v } \leq \sum_{n \geq 0} \theta^n D_\epsi(x)  \abs{v}^{q+1} =  \varphi_\epsi(x) \abs{v}^{q+1}. \]
Thus, by the very definition of
$\rho_\epsi$, the commutative diagram holds and $\TT$ is $\kappa$-tangent to $\II$. The
 commutative relation  $T_x =  N_{\tau(x)} ^{-1} T_{\tau(x)}
F_x$ implies for all $v \in E_x(\rho_\epsi)$ :
\[ \abs { N_x^n T_x(v) - F_x^n(v) } =  \abs {
  T_{\tau^n(x)} F^n_x(v) - F_x^n(v) } \leq \varphi_\epsi(\tau^n(x)) \abs{F^n_x(v)}^{q+1}. \]
By the $\epsi$-fast property of $\varphi_\epsi$, (\ref{ineg}) and the stability of $E(\rho_\epsi)$,
the right hand side is lower than :
\[ \varphi_\epsi(x) e^{n\epsi} (M^n e^{n\epsi}) ^{q+1} \abs{v}^{q+1}  =  \varphi_\epsi(x) (m \beta )^n \abs{v}^{q+1}   \leq \varphi_\epsi(x) (m\theta)^n \abs{v}^{q+1},  \] 
which completes the proof of the theorem \ref{HOC}. \\

We shall now define the functions $r_\epsi$, $D_\epsi$
and establish the estimates (\ref{ineg}) and (\ref{tyu}). Let $\epsi' := {\epsi \over q+2}$. By the lemma
\ref{NN}(1), there exists a $\epsi'$-slow function $\phi_{\epsi'}$ such that
(\ref{ineg}) and the following property  hold on $E(\phi_{\epsi'})$ : 
\begin{equation}\label{ineg2} \forall (u,v) \in E_x(\phi_{\epsi'}) \ , \ m e^{-\epsi} \abs{u-v} \leq \abs{N_x(u) - N_x(v)} \leq M e^{\epsi} \abs{u-v}.
\end{equation}
Define $C_\epsi := 4  / \phi_{\epsi'} ^{q+1}$, $D_\epsi := C_\epsi / m e^{-\epsi}$ and 
\[ r_\epsi := \min \left \{ { \phi_{\epsi'} \over 2 } \ , \  { m e^{-\epsi}\phi_{\epsi'} \over \si C_\epsi} \ , \ { m e^{-\epsi} \phi_{\epsi'} \over M
  e^{\epsi} }   \right \}.\]
The functions $C_\epsi$, $D_\epsi$ are ${q+1 \over q+2} \epsi$-fast
and $r_\epsi$ is $\epsi$-slow. Since $r_\epsi \leq \phi_{\epsi'}$, the estimate (\ref{ineg}) is
satisfied on $E(r_\epsi)$. 
In particular, as $Me^{\epsi}\le e^{-\epsi}$, $E(r_{\epsi})$ is stable by $\FF$ and $\NN$.
Also it follows from (\ref{ineg2}) that $N_x$ is invertible on $E_{\tau(x)}(m  e^{-\epsi} \phi_{\epsi'}(x))$ and satisfies :
\begin{equation}\label{inverse}
 \forall (u',v') \in E_{\tau(x)}(m  e^{-\epsi} \phi_{\epsi'}(x)) \ , \ \abs{N_{\tau(x)}^{-1}(u') - N_{\tau(x)}^{-1}(v')} \leq {1\over m e^{-\epsi}} \abs{u'-v'}.
\end{equation} 
In order to establish (\ref{tyu}) we need the following lemma. The last point relies on the tangency of $\FF$ and $\NN$ at the order $q+1$.

\begin{lem}\label{rond}
For all $v \in E_x(r_\epsi)$ and $n \geq 0$ : 
\begin{enumerate}
\item
$\left[ \sum_{j=0}^n C_\epsi(\tau^j(x)) \beta^j \right]   \abs{v} \leq m e^{-\epsi} \phi_{\epsi'}(x)$.
  
\item
$ \max \left\{  \abs{F_x(v)} \ , \  \abs{ N_x(v) } \right\} \leq  m e^{-\epsi} \phi_{\epsi'}(x)$. 

\item
$ \abs{ F_x(v)  -  N_x(v)  } \leq C_\epsi(x)  \abs{v}^{q+1}$.
\end{enumerate} 
\end{lem}

\proof By the $\epsi$-fast property
of $C_\epsi$, the left hand side in the point 1 is lower than   
\[  \sum_{j=0}^n  (\beta e^\epsi) ^j C_\epsi(x) \abs{v} \leq  \si  C_\epsi(x) r_\epsi(x) \leq  \si C_\epsi(x) { m e^{-\epsi} \phi_{\epsi'}(x) \over \si C_\epsi(x) }  =  m e^{-\epsi} \phi_{\epsi'}(x).  \]
The point 2 is a consequence of (\ref{ineg}) :
\[ \max \left\{  \abs{F_x(v)} \ , \  \abs{ N_x(v) } \right\} \leq  M e^{\epsi} r_\epsi(x) \leq M e^{\epsi} {m e^{-\epsi} \phi_{\epsi'}(x) \over M
  e^{\epsi}}  = m e^{-\epsi}  \phi_{\epsi'}(x).  \]
We now come to the last point. By lemma \ref{cow} and the stability of $E(r_\epsi)$ by $\FF$ and $\NN$, we have for all $j \geq 1$ and $v \in E_x(r_\epsi)$ : 
\[ \max \left\{ \abs{F_x^{(j)}(v)}
\ , \  \abs{N_x^{(j)}(v)}  \right\} \leq  {  \phi_{\epsi'}(\tau(x)) \over \phi_{\epsi'}^j(x) }   \abs{v}^{j} \leq \left(  {  \abs{v} \over  \phi_{\epsi'}(x) } \right)^j \leq  \left({r_\epsi(x) \over \phi_{\epsi'}(x)} \right)^j \leq  {1 \over 2^j} . \]
As $\FF$ and $\NN$ are tangent at the order $q+1$, we deduce that :
\[ \abs{ F_x(v)  -  N_x(v)  } \leq 2   \left( {  \abs{v} \over
  \phi_{\epsi'}(x)} \right)^{q+1}   \sum_{m \geq 0} {1 \over
  2^m} =  4  { \abs{v}^{q+1}  \over \phi^{q+1}_{\epsi'}(x) }  = C_\epsi(x) \abs{v}^{q+1}. \ \ \ \ \ \textrm{\fin} \]

To end the proof of the theorem,  we establish (\ref{tyu}) by induction. Let us rewrite (\ref{tyu}) explicitely :
\[ (p_n) \ : \ \forall v \in E_x(r_\epsi) \ , \ \abs{ N^{-(n+1)}_{\tau^{n+1}(x)} F_x^{n+1} (v) - N_{\tau^n(x)}^{-n} F^n_x (v) } \leq \beta^n D_\epsi(\tau^n(x)) \abs{v}^{q+1}. \]
The proof of $(p_n) \Rightarrow (p_{n+1})$ will use the following auxiliary inequality :
\[ (q_n) \ : \  \forall v \in E_x(r_\epsi) \ , \  \abs{ N^{-n}_{\tau^{n+1}(x)} F^{n+1}_x (v) }
  \leq  \left[ \sum_{j=0}^n C_\epsi(\tau^j(x)) \beta^j  \right]  \abs{v} .  \]
The assertion $(p_0)$ is a consequence of (\ref{inverse}) and the two last points of lemma
\ref{rond} : 
\[  \abs{ N^{-1}_{\tau(x)}
  \left( F_x (v) \right) -  N^{-1}_{\tau(x)} \left( N_x (v)
  \right)}   \leq  {1\over m e^{-\epsi}} \abs{ F_x(v)-N_x(v)
  } \leq D_\epsi(x)  \abs{v}^{q+1} .    \]
The assertion $(q_0)$ follows from (\ref{ineg}) and the observation
  $Me^\epsi \leq 1 \leq C_\epsi$. \\

Assume now that $(p_n)$ and $(q_n)$ are satisfied. Let $v' \in E_{x'}(r_\epsi)$, $v := F_{x'}(v')$ and $x :=
\tau(x')$. Observe
that $v \in E_x(r_\epsi)$ because the tube $E(r_\epsi)$ is stable by
$\FF$. Using $D_\epsi = C_\epsi / m e^{-\epsi}$ and (\ref{ineg}), the assertion
$(p_n)$ yields :
\begin{displaymath}
\begin{array}{rcl}
  \abs{ N^{-(n+1)}_{\tau^{n+2}({x'})} F^{n+2}_{{x'}}({v'}) -
  N^{-n}_{\tau^{n+1}({x'})} F^{n+1}_{{x'}} ({v'}) } & \leq  &  \beta^n {C_\epsi(\tau^{n+1}({x'})) \over m e^{-\epsi}}
 \abs{F_{x'}(v')}^{q+1} \\
                                                    & \leq  &  \beta^n {C_\epsi(\tau^{n+1}({x'})) \over m e^{-\epsi}}
 (M e^\epsi)^{q+1} \abs{{v'}}^{q+1},
\end{array}                                      
\end{displaymath}
which, by the definition of $\beta$, leads to :
\begin{equation}\label{interm}
  \abs{ N^{-(n+1)}_{\tau^{n+2}({x'})} F^{n+2}_{{x'}}({v'}) -
  N^{-n}_{\tau^{n+1}({x'})} F^{n+1}_{{x'}} ({v'}) }  \leq  \beta^{n+1}  C_\epsi(\tau^{n+1}({x'})) \abs{{v'}}^{q+1}.
\end{equation}
We now deduce  $(q_{n+1})$ from $(q_n)$ and (\ref{interm}). For all
$v' \in E_{x'}(r_\epsi)$ we have :
\begin{displaymath}
\begin{array}{rcl}
    \abs{ N^{-(n+1)}_{\tau^{n+2}({x'})}
  F^{n+2}_{x'}({v'})}   & \leq  &   \beta^{n+1} C_\epsi(\tau^{n+1}({x'})) \abs{{v'}} ^{q+1} +  \left[
  \sum_{j=0}^n C_\epsi(\tau^j({x'})) \beta^j  \right]  \abs{{v'}}  \\
       &  \leq &  \left[
  \sum_{j=0}^{n+1} C_\epsi(\tau^j({x'})) \beta^j  \right]  \abs{{v'}}.
\end{array}                                      
\end{displaymath}
It follows from $(q_n)$, $(q_{n+1})$ and the lemma \ref{rond}-(1) that 
$\abs{ N^{-(n+1)}_{\tau^{n+2}({x'})} F^{n+2}_{{x'}}({v'}) }$ and
$\abs{ N^{-n}_{\tau^{n+1}({x'})} F^{n+1}_{{x'}} ({v'}) }$ are lower
than $m e^{-\epsi} \phi_{\epsi'}(x')$. We may thus compose (\ref{interm}) by
$N^{-1}_{\tau({x'})}$ and, using (\ref{inverse}), get :
\begin{displaymath}
\begin{array}{rcl}
 \abs{ N^{-(n+2)}_{\tau^{n+2}({x'})} F^{n+2}_{{x'}}({v'}) -
                                                          N^{-(n+1)}_{\tau^{n+1}({x'})}
                                                          F^{n+1}_{{x'}}
                                                          ({v'})}  &
                                                          \leq & \beta^{n+1}  {C_\epsi(\tau^{n+1}({x'})) \over m e^{-\epsi}} \abs{{v'}}^{q+1} \\
                                                          & \leq & \beta^{n+1}  D_\epsi(\tau^{n+1}({x'})) \abs{{v'}}^{q+1},
\end{array}                                      
\end{displaymath}
which is the assertion $(p_{n+1})$. \fin

\subsection{Normalization at a finite order}\label{fo}

We prove here that a bundle map which satisfies the assumptions of theorem \ref{normaliza} is, for any prescribed order, tangent to some resonant 
bundle map :
 
\begin{prop}\label{poincare}
Let $\GG$ be a tame and regular contracting bundle map over
$\tau$. For any $p \geq 1$, there exists a resonant bundle map $\RR$
such that $\GG$ is tangent to $\RR$ at the order $p+1$. More precisely,
for any $\epsilon >0$ and $0<\kappa<1$, there exist a $\epsi$-slow function $\si_\epsi$, a
bundle map $\UU$ which is $\kappa$-tangent to $\II$ and a bundle map $\FF = \RR +
O(p+1)$ such that the following diagram 
commutes :
\[
\xymatrix{
          E(\si_\epsi) \ar[r]^\GG \ar[d]_\UU   & E(\si_\epsi) \ar[d]^\UU \\
           E(2\si_\epsi)  \ar[r]^\FF &      E(2\si_\epsi)         }
\]
 and the tubes $E(\si_\epsi), E(2\si_\epsi)$ are contracted.
\end{prop}

Following Jonsson-Varolin (\cite{JV}, lemma 5.2), we first solve the Poincaré's homological equation in a non-autonomous setting. 

\begin{prop}\label{coho}
Let $\AA$ be a linear regular contracting bundle map (see
definition \ref{rc}). Let $\HH = \HH^{(m)}$ be a $m$-homogeneous bundle map over $\tau$ and
$\RRR(\HH)$ be its resonant part. Assume that $\abs \HH$ is a fast
function. Then there exists a $m$-homogeneous bundle map $\QQ$ over $\Id_X$ such that $\abs
\QQ$ is a fast function and $\HH + ( \QQ \circ \AA - \AA \circ
\QQ) = \RRR(\HH)$.
\end{prop}

\proof
It suffices to consider a bundle map $\MM
:= \pi_j^\al(\HH)$, where $\abs{\al} = m$ and $1 \leq j \leq k$. We treat the three cases $\al \in \RRR_j,  \mathfrak{B}_j$ and
$\mathfrak{P}_j$ separately. If $\al \in \RRR_j$, we let $\QQ = 0$. If $\al \in  \mathfrak{B}_j$, we set $\QQ  :=  \sum_{n \geq 0} \AA^{-(n+1)}
\MM \AA^n$, i.e. 
\begin{equation}\label{sss}  
\forall x \in X \ , \ Q_x =  A^{-1}_{\tau(x)} M_x + \sum_{n \geq 1} \left[  A^{-1}_{\tau(x)} \ldots A^{-1}_{\tau^{n+1}(x)} \,
M_{\tau^{n}(x)} \, A_{\tau^{n-1}(x)}  \ldots A_{x}\right] .
\end{equation}
A formal computation shows that $\QQ$ is $m$-homogeneous and satisfies $M_x + (Q_{\tau(x)} \circ A_x - A_x \circ
Q_x  ) = 0 = \RRR ( M_x )$ for all $x \in X$. We now check the convergence of the series
(\ref{sss}), which we denote $Q_x :=  \sum_{n \geq  0} Q_{n,x}$. Let
$\psi_\epsi  \geq \abs{\HH} \geq \abs{\MM}$ be a $\epsi$-fast function. Observe that each of the $k_j$ coordinates of $\MM : \XX \to \LL^j$ is a 
multiple of the monomial $P_\al$. We have therefore by the
$\epsi$-fast property of $\psi_\epsi$ and the regular contracting
property of $\AA$ :
\begin{displaymath}
\begin{array}{rcl}
\forall n \geq 1 \ , \   \abs{ Q_{n,x}(v) }     &  \leq   &   e^{-(n+1)(\Lambda_j - 
  \epsilon)}  \,  \left[\psi_\epsi(x)  e^{n\epsilon} \right] \,
e^{n(\lambda \cdot \al + \abs \al \epsilon) } \,  \abs{v}^m  \\
       &  \leq   &  e^{n [ \lambda \cdot\al -
  \Lambda_j +  (\abs \al + 2) \epsilon ] } e^{-\Lambda_j + \epsi} \psi_\epsi(x) \abs{v}^m \\
  &  \leq   &   e^{-n \zeta}  e^{-\Lambda_j + \epsi} \psi_\epsi(x) \abs{v}^m,
\end{array}                                      
\end{displaymath} 
where $\zeta >0$ is defined in the lemma \ref{marge}. Besides the
convergence of (\ref{sss}) this shows that $\abs \QQ \leq  [ e^{-\Lambda_j +
  \epsi}  \sum_{n \geq 0} e^{-n \zeta} ] \psi_\epsi(x)$, so that $\abs \QQ$
is fast. 

Finally, if $\al \in \mathfrak{P}_j$, we set $\QQ  :=  - \sum_{n
  \geq 0}  \AA^n \MM \AA^{-(n+1)}$, i.e. :
\[  Q_x = - M_{\tau^{-1}(x)} A^{-1}_x - \sum_{n \geq 1} \left[ A_{\tau^{-1}(x)}  \ldots A_{\tau^{-n}(x)}  \,
M_{\tau^{-(n+1)}(x)} \,  A^{-1}_{\tau^{-n} (x)} \ldots
  A^{-1}_x \right] . \]
We obtain in that case :
\begin{displaymath}
\begin{array}{rcl}
   \abs{ Q_{n,x}(v) }    &  \leq   &   e^{n(\Lambda_j + 
  \epsilon)}  \,  \left[\psi_\epsi(x)  e^{(n+1)\epsilon} \right] \,
e^{(n+1)(-\lambda \cdot \al + \abs \al \epsilon) } \,  \abs{v}^m  \\
       &  \leq   &  e^{(n+1) [ -\lambda \cdot\al +
  \Lambda_j +  (\abs \al + 2) \epsilon ] } e^{-(\Lambda_j + \epsi)} \psi_\epsi(x) \abs{v}^m \\
  &  \leq   &   e^{-(n+1) \zeta}  e^{-(\Lambda_j + \epsi)} \psi_\epsi(x) \abs{v}^m ,
\end{array}                                      
\end{displaymath} 
and we conclude as before. \fin

\begin{cor}\label{etapem}
Let  $\KK : E(s_\epsi) \to E(s_\epsi)$ be a regular contracting
bundle map, where $s_\epsi$ is $\epsi$-slow. Suppose that $\Lip(\KK) \leq c < 1$ on $E(s_\epsi)$. Let $m \geq 2$ and $\gamma > 0$ small 
(depending on $c ,\epsi$). There exist a $\epsi$-slow
function $\tilde{r_\epsi} \leq s_\epsi$, a bundle map $\SS_m$ tangent to $\II$ and a bundle map ${\widehat \KK}$
tangent to $\AA + \KK^{(2)}  + \ldots +  \KK^{(m-1)} + \RRR(\KK^{(m)})$ at the order $m+1$, 
such that the following diagram commutes for any $\epsi$-slow function $r_\epsi \leq \tilde{r_\epsi}$ :
\[
\xymatrix{
          E(e^{-\gamma} r_\epsi) \ar[r]^\KK \ar[d]_{\SS_m} &
          E(e^{-\gamma} r_\epsi) \ar[d]^{\SS_m} \\
           E(r_\epsi)  \ar[r]^{\widehat \KK}     &  E(r_\epsi)\;.         }
\]
 Moreover $\Lip (\widehat{K}) \leq c e^{2\gamma}$. 
 \end{cor}

\proof 
 Let $\gamma >0$ such that $c e^{2\gamma} \leq
e^{-\epsi}$. 
The bundle map $\KK$ is tame, so the function $\abs{\KK^{(m)}}$ is fast (see lemma
\ref{cow}) and we may apply the proposition \ref{coho} : there exists a $m$-homogeneous bundle map $\QQ$
over $\Id_X$ such that $\abs{\QQ}$ is fast and 
\begin{equation}\label{coh}
  K_x^{(m)} + ( Q_{\tau(x)} \circ A_x - A_x \circ Q_x) = \RRR (K^{(m)}_x).
\end{equation}
We let $\SS := \II + \QQ$, this bundle map is tame. As $\SS$ is tangent to $\II$, one may decrease $\tilde{r_\epsi}$ to obtain (see lemma \ref{NN}(1)) : 
\[ \forall (u,v) \in E_x( e^{\gamma} \tilde{r_\epsi}) \ , \ e^{-\gamma} \abs{u-v} \leq
\abs{S_x(u) - S_x(v)} \leq e^\gamma \abs{u-v}. \] 
In particular, $\SS^{-1}$ exists on $E(\tilde{r_\epsi})$. Thus, for any $\epsi$-slow function
$r_\epsi \leq \tilde{r_\epsi}$, the bundle map $\widehat \KK := \SS \circ \KK \circ \SS^{-1}$ is well defined on
$E(r_\epsi)$ and is $c e^{2\gamma}$-Lipschitz. Furthermore, $E(r_\epsi)$ is
stable by $\widehat \KK$ since $c e^{2\gamma} \leq e^{-\epsi}$. It remains to prove that $\widehat \KK$ is tangent to $\AA + \KK^{(2)}  + \ldots +
\KK^{(m-1)} + \RRR (\KK^{(m)})$ at the order $m+1$. Let us write $U_x
\simeq V_x$ when $U_x - V_x = O(m+1)$. Observe first that $S_x^{-1}
\simeq \Id_{\C^k} - Q_x$. Moreover, as $\QQ$ is homogeneous of degree $m \geq
2$, we have :
\begin{displaymath}
\begin{array}{rcl}
      K_x S_x^{-1}  &  \simeq  &  \left[ ( A_x +  K_x^{(2)} + \ldots +
K_x^{(m-1)} ) + K_x^{(m)} \right] \circ \left(  \Id_{\C^k}  - Q_x  \right) \\
                   &  \simeq  &  ( A_x +  K_x^{(2)} + \ldots +
K_x^{(m-1)} ) + K_x^{(m)} - A_x \circ Q_x
\end{array}                                      
\end{displaymath} 
and then :
\begin{displaymath}
\begin{array}{rcl}
   S_{\tau(x)} K_x S_x^{-1}  &   \simeq  &\left(  \Id_{\C^k}+
   Q_{\tau(x)} \right)
\circ  \left[   ( A_x +  K_x^{(2)} + \ldots +  K_x^{(m-1)} ) + K_x^{(m)} - A_x \circ Q_x     \right]  \\
                   &   \simeq     &  ( A_x +  K_x^{(2)} + \ldots
+  K_x^{(m-1)} ) + K_x^{(m)}  +  Q_{\tau(x)} \circ  A_x - A_x \circ
Q_x  \\
              &   \simeq     &  ( A_x +  K_x^{(2)} + \ldots +  K_x^{(m-1)} ) + \RRR (K_x^{(m)}),
\end{array}                                      
\end{displaymath} 
where the last line follows from (\ref{coh}).   \fin

\proo {\sc of the proposition \ref{poincare}} : If $p = 1$, we just take
$\RR = \AA$, so let $p \geq 2$. By the lemma \ref{NN}(1), there exists
a $\epsi$-slow function $s_\epsi$ such that $\GG : E(s_\epsi) \to
E(s_\epsi)$ and $\Lip {\GG} \leq e^{\Lambda_1+\epsilon} e^\epsi < 1$ (take $b = e^{\Lambda_1 + \epsi}$ and $\ga = \epsi$). Let  $c' := e^{\Lambda_1+2\epsilon}$ and $\ga$ small enough
such that $1/2 \leq e^{-(p-1)\gamma}$ and $c' e^{2({p}-1)\gamma} < 1$. We apply successively the
corollary \ref{etapem} with $(m,c) = (2,c'),(3,c'e^{2\gamma}),\ldots,(p,c'e^{2(p-2)\gamma})$. We obtain bundle maps $\SS_2,\SS_3,\ldots,\SS_p$ tangent to $\II$ and a $\epsi$-slow function $r_\epsi$ such that the following diagram is commutative, where $\UU := \SS_{p} \circ
\ldots \circ \SS_2$ :
\[
\xymatrix{
          E(e^{-(p-1)\gamma} r_\epsi) \ar[r]^\GG \ar[d]_{\UU} &
          E(e^{-(p-1)\gamma} r_\epsi) \ar[d]^{\UU} \\
           E(r_\epsi)  \ar[r]^{\FF}     &  E(r_\epsi)         }
\]
The bundle map $\FF$ is tangent to $\RR$ at the order $p+1$ and
satisfies $\Lip(\FF) \leq c' e^{2({p}-1)\gamma} < 1$ on
$E(r_\epsi)$. We have $\UU : E(r_\epsi/2) \to E(r_\epsi)$ because $1/2
\leq e^{-(p-1)\gamma}$. We set $\si_\epsi := r_\epsi /2$. The bundle
map $\UU$ is $\kappa$-tangent to $\II$ on $E(\si_\epsi)$ (decrease $\si_\epsi$ if necessary, see lemma \ref{NN}(2)). \fin 

\subsection{Proof of theorem \ref{normaliza}}

 Let $\GG$ be a tame and regular contracting bundle map over $\tau$, let $m := e^{\Lambda_l - \epsi}, M := e^{\Lambda_1 +
 \epsi}$ and ${\tilde q}$ be the entire part of $\Lambda_l /
 \Lambda_1$ (so that $M ^{{\tilde q}+1} < m$). We apply successively
 the proposition \ref{poincare} and the theorem \ref{HOC} with $p=q = {\tilde
 q}$. There exists bundle maps $\TT,\UU$ tangent to $\II$, a resonant
 bundle map $\RR$ and a $\epsi$-slow function $\si_\epsi$  such that
 the following diagram  commutes :
\[
\xymatrix{
          E(\si_\epsi)   \ar[rr]^\GG              \ar[d]_\UU      &        &   E(\si_\epsi) \ar[d]^\UU \\
          E(2\si_\epsi)  \ar[rr]^{\RR +  O({{\tilde q}}+1)} \ar[d]_\TT     &        &   E(2\si_\epsi) \ar[d]^\TT  \\
          E(4\si_\epsi)  \ar[rr]^\RR                         &         &   E(4\si_\epsi)                   }
\]
We let $\VV := \TT \circ \UU$ and $\rho_\epsi \leq \si_\epsi$ be a $\epsi$-slow function such that $\VV$ is 
$\kappa$-tangent to $\II$ on $E(\rho_\epsi)$ (see lemma \ref{NN}). We obtain $\VV : E(\rho_\epsi) \to E(2\rho_\epsi)$ by choosing $\kappa < 1$.\fin

\section{Application to holomorphic dynamics}

We prove in this section the theorems \ref{nt} and \ref{appli}, let us first recall the setting in which we will work.
We consider an  
holomorphic endomorphism $f : \Pj^k \to \Pj^k$ of algebraic degree $d \geq 2$. The topological
degree of $f$ is $d_t := d^k$. 
The equilibrium measure $\mu$ of $f$ is given by $\mu = \lim_{n \to \infty} {1 \over d_t^n} {f^n}^* \om^k$, where $\om$ is the Fubini-Study
$(1,1)$ form on $\Pj^k$. This measure is mixing, satisfies $\mu(f(B))
= d_t \mu(B)$ whenever $f$ is injective on the borel set $B$, and does
not charge any analytic subset of $\Pj^k$. As Briend-Duval \cite{BD}
proved, the Lyapunov exponents $\chi_1 \leq \ldots \leq \chi_k$
of $\mu$ are bounded from below by $\log \sqrt d$. We are interested
in the following quantities : $$\Sigma_s :=
\chi_{k-s+1} + \ldots + \chi_k.$$ 

Let $\orb := \{ \hat x := (x_n)_{n
  \in \Z} \, , \,  x_{n+1} = f (x_n) \}$ be the set of orbits. We note $\hat f$ the left shift sending
$(\ldots,x_{-1},x_0,x_1,\ldots)$ to
$(\ldots,x_{0},x_1,x_2,\ldots)$ and $\tau := \hat f^{-1}$. Let $\pi$ be the projection $\hat x \mapsto x_0$ and $\nu$ be the unique probability measure on $\orb$ 
invariant by $\tau$ and satisfying $\nu (\pi ^{-1}
(B)) = \mu(B)$ on the borel sets $B \subset \Pj^k$. The measure $\nu$ is mixing. We will work with the
subset $X := \{ \hat x \in  \orb  \ , \  x_n \notin \CC_f \ , \
\forall n \in \Z \}$ where $\CC_f$ is the critical set of $f$. This subset has full $\nu$-mesure, because $\mu(\CC_f)
= 0$. For all $\hat x \in X$, the inverse
branch of $f^n$ that sends $x_0$ to $x_{-n}$ is denoted $f^{-n}_{\hat x}$. \\

We note $d(.,.)$ the distance on $\Pj^k$ induced by the Fubini-Study metric, $B_x(r) \subset \Pj^k$ the ball centered at $x$ of radius $r$ for this metric 
and $\B(r)$ the ball $\{ \abs{z} < r \} \subset \C^k$ for the standard metric.

\subsection{Normalization along orbits of endomorphisms of $\Pj^k$} \label{nts}

Our aim here is to prove the theorem \ref{nt}. To this purpose we will first construct the bundle map $\FF^{-1}$ of the inverse branches of an  
endomorphism $f : \Pj^k \to \Pj^k$. We will then apply to $\FF^{-1}$ the Oseledec-Pesin reduction theorem and our theorem \ref{normaliza}.
For defining this bundle map we shall use 
a family of charts $(\psi_x)_{x \in \Pj^k}$ satisfying the following properties, where $0 < M_0 < 1$ is a constant independant of $x \in \Pj^k$ :

{(P1) } $\psi_x : \C^k \to \Pj^k$ is a biholomorphism onto its image and $\psi_x(0)=x$.

{(P2) } for all $(z_1,z_2) \in \B(M_0)$, ${1 \over 2} \abs{z_1-z_2} \leq d(\psi_x(z_1) , \psi_x(z_2)   ) \leq 2 \abs{z_1-z_2}$.

\noindent We also require the following condition where 
$\abs{\,.\,}_{x,u}$ is the norm on the spaces
$(\bigwedge^s \C^k)_{1 \leq s \leq k}$ inherited by the strictly
positive $(1,1)$ form $(\psi_x^*\om)(u)$ :

{(P3) } for all $u \in \B(M_0)$, ${1 \over 10}  \abs{\,.\,} \leq \abs{\,.\,}_{x,u}
  \leq 10 \abs{\,.\,}$ and $\abs{\,.\,}_{x,0} =
  \abs{\,.\,}$. \\

For $\hat x \in X$, we identify $E_{\hat x} =  \{(\hat x,v) \
, \ v \in \C^k \}$ with $\C^k$ and let $\psi_{\hat x} : E_{\hat x} \to \Pj^k$ be the map $\psi_{\hat x} := \psi_{x_0} : \C^k \to \Pj^k$. We define
\[ F_{\hat
 x} := \psi^{-1}_{{\hat f(\hat x)}} \circ f \circ \psi_{{\hat x}}. \]
This map satisfies  $F_{\hat
 x}(0) = 0$. By the uniform continuity of $f$, there exists $0 < M_1 \leq M_0$ such that 
the following bundle map is well defined :
\begin{displaymath} 
 \FF :
\begin{array}{rccc}
 &  E(M_1)          &  \longrightarrow    &         E(M_0)   \\
 &  (\hat x , z )   & \longmapsto &  \left(\hat f(\hat x), F_{\hat
 x}(z)  \right).
\end{array}                                      
\end{displaymath}

Recall that $X = \{ (x_n)_{n \in \Z} \, , \,  x_{n+1} = f (x_n) \, , \,  x_n \notin \CC_f \}$. So for any $\hat x \in X$ the map $F_{\tau(\hat x)}$ 
is invertible in a neighbourhood of the origin. As the next lemma shows, the bundle map $\FF^{-1}$ is actually tame.

\begin{lem}\label{abd}
There exists a slow function $\al : X \to ]0,1]$ and a (tame) bundle map
\begin{displaymath} 
 \FF^{-1} :
\begin{array}{rccc}
 &  E(\al)          &  \longrightarrow    &         E(M_1)   \\
 &  (\hat x , z )   & \longmapsto &  \left( \tau(\hat x), F^{-1}_{\hat x}(z) \right),
\end{array}                                      
\end{displaymath}  
where $F^{-1}_{\hat x}:=  
( F_{\tau(\hat x)})^{-1} = \left( \psi^{-1}_{\hat x} \circ f \circ \psi_{\tau(\hat x)} \right)^{-1}$.  
\end{lem}

\proof Let $t(\hat x) := \norm {  ( d_{x_{-1}} f  )^{-1} }^{-2}$. There
exists a constant $c >0$ depending only on the first and second
derivatives of $f$ such that the map $F^{-1}_{\hat x}$ exists on
$E_{\hat x}(c t)$. We use here a quantitative version of the
inverse local theorem (see e.g. \cite{BD}, lemma 2). We let $\al :=
\min \{
c t, 1 \}$. As $\log \al$ is $\nu$-integrable (see \cite{S},
subsection 3.7), there exists a $\epsi$-slow function
$\al_\epsi : X \to ]0,1]$ such that $\al_\epsi \leq \al$ (see lemma \ref{abelard}). The function
$\al$ is therefore slow. \fin

Let $\DD$ be the linear part of $\FF^{-1}$. The Oseledec-Pesin theorem asserts that $\DD$ is regular contracting
after conjugation by a ``tempered'' family $\CC$ of linear maps :
 
\begin{thm}[Oseledec-Pesin $\epsi$-reduction \cite{KH}]\label{OPR}
There exist a linear bundle map $\CC$ over $\Id_X$ and a function $h_\epsi : X \to [1,+\infty[$ such that :
\begin{enumerate}
\item $\AA  := \CC \circ \DD \circ \CC^{-1}$ is regular contracting.
\item $\forall \hat x \in X, \, \abs{v} \leq
\abs {C_{\hat x} (v)}  \leq  h_\epsi(\hat x) \abs{v}$ and 
$e^{-\epsi} h_\epsi(\hat x) \leq h_\epsi(\tau (\hat x)) \leq e^{\epsi}
h_\epsi(\hat x)$.
\end{enumerate}
\end{thm}
The contraction rates $\Lambda_l < \ldots < \Lambda_1 <0$ of $\AA$ are in our context the distinct
opposite Lyapunov exponents of $(\Pj^k,f,\mu)$ (i.e. the distinct $-\chi_i$) and are negative. The integers $k_j$ are the multiplicities of these exponents. \\

As the reader may easily check, the following proposition is a version of
the theorem \ref{nt}. We stress that the points 2 and 3 are consequences of the
algebraic properties of the resonant maps (see proposition
\ref{estreso}). We recall that $\epsi \ll \abs{\Lambda_1}$.

\begin{prop}\label{propbd}
With the preceding notations, there exist  $\epsi$-slow
functions $\eta_\epsi , t_\epsi : X \to ]0,M_0]$, a resonant bundle map $\RR$ over $\tau$ and an injective bundle map $\WW$ over $\Id_X$ (tangent to $\CC$) such that the following
diagram commutes for all $n \geq 1$ :
\[
\xymatrix{
          E(\eta_\epsi) \ar[r]^{\FF^{-n}} \ar[d]_\WW   & \FF^{-n} E(\eta_\epsi)   \ar[d]^\WW \\
           E(t_\epsi)  \ar[r]^{\RR^n}         &  E(t_\epsi)\;\;\;.        }
\]
We have $\FF^{-n} E(\eta_\epsi) \subset E(M_0)$. There exist $\al' >0$ and $\epsi$-fast functions $\beta'_\epsi, L_\epsi',T_\epsi' : X \to [1,+\infty[$ such that :
\begin{enumerate}
\item $\forall (y,y') \in  F_{\hat x}^{-n}(E_{\hat x}(\eta_\epsi))$, $\al' \, \abs{y-y'} \leq \abs{W_{\tau^n(\hat x)}(y) - W_{\tau^n(\hat x)}(y') } \leq  \beta'_{\epsi}(\tau^n(\hat x)) \abs{y-y'}$.
\item $\Lip F_{\hat x}^{-n} \leq L_\epsi'(\hat x) e^{-n \chi_1 + n \epsi}$
  on $E_{\hat x}(\eta_\epsi)$,
\item for all $z \in F_{\hat x}^{-n}(E_{\hat x}(\eta_\epsi)), \, \vert
  {1 \over n} \log \abs {\bigwedge^s d_z F_{\tau^n(\hat x)}^n}   -
  \Sigma_s  \vert  \leq {1 \over n} \log T_\epsi'(\hat x) + \epsilon$.
\end{enumerate}
\end{prop}

\proof Let $\epsi' = \epsi
/2(k+\eta)$ (the constant $\eta$ is defined in the proposition \ref{estreso}
and depends only on ${\tilde q},k$). By the lemma \ref{abd} and the
theorem \ref{OPR}, there exist a ${\epsi'}$-slow function $\al_{\epsi'} \leq \al$  and a ${\epsi'}$-fast function $h_{\epsi'}$ such that the bundle map $\GG := \CC \circ \FF^{-1} \circ \CC^{-1}$ is well defined $E(\al_{\epsi'}) \to E(M_1 h_{\epsi'})$ and is regular contracting (observe that $E(\rho) \subset \CC(E(\rho)) \subset E(\rho. h_{\epsi'})$). In particular $\GG$ is tame. By the theorem
\ref{normaliza}, there exists a ${\epsi'}$-slow function
$\xi_{\epsi'}$, a bundle map $\VV$ $1/2$-tangent to $\II$ and a resonant bundle map $\RR$ such the following diagram commutes :
\[
\xymatrix{
          E(\eta_{\epsi})   \ar[rr]^{\FF^{-n}} \ar[d]_{\CC}      &
          &  \FF^{-n}  E(\eta_{\epsi})  \ar[d]^{\CC}  \\
          E(\xi_{\epsi'})  \ar[rr]^{\GG^n}           \ar[d]_\VV     &        &   E(\xi_{\epsi'}) \ar[d]^\VV  \\
          E(2\xi_{\epsi'})  \ar[rr]^{\RR^n}                         &         &   E(2\xi_{\epsi'})\;\;\; .                   }
\]
We may assume that $\xi_{\epsi'} \leq M_0/2$ and that $\abs{
  \bigwedge^s \Id_{\C^k} - \bigwedge^s (d_t V_{\hat x})^{\pm 1} } \leq
1/10$ for any $t \in E_{\hat x}(\xi_{\epsi'})$ and $\hat x \in X$ (lemma \ref{NN}(3)). Let $\eta_{\epsi} = \xi_{\epsi'} / h_{\epsi'}$ and $t_{\epsi} := 2
\xi_{\epsi'}$. These functions are $\epsi$-slow and satisfy $\eta_{\epsi} , t_{\epsi} \leq M_0$. Moreover $\FF^{-n}  E(\eta_{\epsi}) \subset \CC^{-1} E(\xi_{\epsi'}) \subset   E(\xi_{\epsi'}) \subset E(M_0)$. We let $\WW := \VV \circ \CC$. Observe that :
\[ \forall (y,y') \in  F_{\hat x}^{-n}(E_{\hat x}(\eta_\epsi)) \ , \ {1 \over 2} \abs{y-y'} \leq \abs { W_{\tau^n(\hat x)}(y) - W_{\tau^n(\hat x)}(y')  } \leq  {3 \over 2} h_{\epsi'}(\tau^n(\hat x)) \abs{y-y'}.  \]
The point 1 follows with $\al' = 1/2$ and $\beta'_\epsi =  {3 \over 2} h_{\epsi'}$. Now we prove the points 2 and 3. Let $y \in E_{\hat x}(\eta_{\epsi})$ and
$t = C_{\hat x}(y)$. We obtain by the theorem
\ref{OPR}(2) and the commutative diagram above : 
\[ \abs{ \bigwedge^s d_{y} F^{-n}_{\hat x} }   \leq
\abs{C_{\tau^n(\hat x)}^{-1} }^s  \abs{ \bigwedge^s d_t G^n_{\hat x} }
\abs{C_{\hat x}}^s \leq \abs{ \bigwedge^s d_t G^n_{\hat x} }
h_{\epsi'}(\hat x)^{s}, \]

\[ \abs{ \bigwedge^s  d_t G^n_{\hat x} }   \leq
\abs{C_{\tau^n(\hat x)}}^s \abs{ \bigwedge^s d_{y} F^{-n}_{\hat x}  }
\abs{C_{\hat x}^{-1}}^s \leq (h_{\epsi'}(\hat x)e^{n{\epsi'}})^s \abs{
  \bigwedge^s d_{y} F^{-n}_{\hat x}}. \]
We deduce the estimate :
\begin{equation}\label{aqw}
 \Big \vert   {1 \over n} \log \abs{ \bigwedge^s d_{y} F^{-n}_{\hat x} }
- {1 \over n} \log \abs{ \bigwedge^s d_t G^n_{\hat x} }  \Big \vert  \leq {1
  \over n} \log  h_{\epsi'}^{k}(\hat x) + k{\epsi'}.
\end{equation}
Now let $w := V_{\hat x}(t)$, $G^n := G^n_{\hat x}$, $R^n :=
R^n_{\hat x}$, $V := V_{\hat x}$, $V_n := V_{\tau^n({\hat x})}$ and write :
\[ \bigwedge^s d_t G^n =  \bigwedge^s ( d_{G^n(t)} V_n )^{-1}   \bigwedge^s d_{w} R^n   \bigwedge^s d_t V  := (\bigwedge^s \Id_{\C^k} + \Omega_1 )    \bigwedge^s d_{w} R^n      ( \bigwedge^s \Id_{\C^k}  +  \Omega_2  ) ,  \] 
where $\Omega_i : \bigwedge^s \C^k \to \bigwedge^s \C^k$ satisfy
  $\abs{\Omega_i} \leq 1/10$. This implies that 
\begin{equation}\label{cvn}
 {1 \over
  2} \abs{ \bigwedge^s d_{w} R^n } \leq \abs {\bigwedge^s d_t G^n }
  \leq  2 \abs{ \bigwedge^s d_{w} R^n }.
\end{equation}
The proposition \ref{estreso} gives a $\eta\epsi'$-fast function $H_{\eta\epsi'}$ such that for all $s \in \{1,\ldots,k \}$ :
\[ \Big \vert   {1 \over n} \log \abs{\bigwedge^s d_w R^n_x} -
(\lambda_1 + \ldots + \lambda_s)  \Big \vert  \leq {1 \over n} \log
H_{\eta \epsi'}(x) + \eta \epsilon'.\]
We deduce with (\ref{aqw}), (\ref{cvn}) and $\lambda_i = -\chi_i$ : 
\begin{equation}\label{vbnn}
 \Big \vert   {1 \over n} \log \abs{ \bigwedge^s d_{y} F^{-n}_{\hat x} }
+ (\chi_1 + \ldots + \chi_s)  \Big \vert  \leq {1
  \over n} \log  [2 H_{\eta \epsi'} h_{\epsi'}^{k} ] (x)  + (k+\eta){\epsi'}.
\end{equation}
The function $L_\epsi'
  := 2H_{\eta \epsi'} h_{\epsi'}^{k}$ is $(k+\eta)\epsi'$-fast (therefore $\epsi$-fast). The point 2
  follows from  (\ref{vbnn}) with $s=1$. For the point 3, observe that $\abs{
  \bigwedge^s d_z F_{\tau^n(\hat x)}^n } =   \abs{ \bigwedge^{k-s} d_y
  F^{-n}_{\hat x} } \abs{\bigwedge^k d_y F^{-n}_{\hat x}}^{-1}$ where
  $z =  F^{-n}_{\hat x}(y)$. We obtain by using (\ref{vbnn}) twice :
\begin{equation*}
 \Big \vert   {1 \over n} \log \abs{ \bigwedge^{k-s} d_{y} F^{-n}_{\hat x} }
+ (\chi_1 + \ldots + \chi_{k-s})  \Big \vert  \leq {1
  \over n} \log L_\epsi' (x)  + (k+\eta){\epsi'},
\end{equation*}
\begin{equation*}
 \Big \vert  {1 \over n} \log \abs{ \bigwedge^k d_{y} F^{-n}_{\hat x} }
+ (\chi_1 + \ldots + \chi_k)  \Big \vert  \leq {1
  \over n} \log L_\epsi' (x)  + (k+\eta){\epsi'}.
\end{equation*}
These estimates imply :
\begin{equation*}
 \Big \vert   {1 \over n} \log \abs{ \bigwedge^s  d_z F_{\tau^n(\hat x)}^n }
- (\chi_{k-s+1} + \ldots + \chi_k)  \Big \vert  \leq {1
  \over n} \log L_\epsi'^2(x)  + 2(k+\eta){\epsi'}.
\end{equation*}
We finally let $T_\epsi' := L_\epsi'^2$, which is a $\epsi$-slow function.\fin

\subsection{An approximation formula for sums of Lyapunov exponents} \label{app}
This subsection is devoted to the proof of theorem \ref{appli}.
Let $\Rep_n$ (resp. $\Rep_n^*$) be the set of repulsive periodic points whose period divides $n$
(resp. equals $n$).
Let us consider the function $\varphi_n$ defined on $\Pj^k$ by
$$\varphi_n (z) := {1 \over n} \log \norm {\bigwedge^s d_z f^n}.$$
Then, the theorem \ref{appli} may be stated as follows:

\begin{equation}\label{exa}
\lim_{n \to
    + \infty}  {1 \over d_t^n} \sum_{p \in \Rep_n} \varphi_n (z)=
\lim_{n \to+ \infty}  {1 \over d_t^n} \sum_{p \in \Rep_n^*} \varphi_n (z) = \Sigma_s.
\end{equation}

Let $\Fix_n := \{ z \in \Pj^k \, , \, f^n(z) = z \}$. Since the number of fixed points of $f^n$
\emph{counted with multiplicity} is, by Bezout's theorem, equal to $1 + d^n + \ldots +
d^{nk}$ (see \cite{S}, subsection
1.3) we have $\cd(\Fix_n) \leq  (k+1)d_t^n$. This implies that for any $n \geq
1$, $\cd(\Rep_n) - \cd(\Rep_n^*) \leq \sum_m \cd(\Fix_m) \leq n (k+1) d_t^{n/2}$, where the sum
runs over the integers $1 \leq m \leq n/2$ which divide $n$. The first equality in (\ref{exa}) is then a consequence of the following lemma.

\begin{lem}\label{l2}
There exists $\Gamma \geq 0$ such that for all $p \in \Rep_n$ and $n
\geq 1$, $0 \leq \varphi_n (p) \leq \Gamma$.
\end{lem}

\proof
Let $p \in \Rep_n$ and $(p_0,p_1,\ldots,p_{n-1})$ be the repulsive
cycle generated by $p_0 := p$. 
Let $\Gamma := k . \max_{z \in \Pj^k} \log^+ \norm {d_z f}$. 
The observations $\norm {\bigwedge^s d_p f^n}
\leq  \norm {d_p f^n}^k  \leq \prod_{i=0}^{n-1}  \norm {d_{p_i} f}^k$ and $0 \leq \varphi_n (p) =  {1 \over n} \log^+ \norm {\bigwedge^s d_p f^n}$ imply the inequalities :
$ 0 \leq \varphi_n (p)  \leq k. {1 \over n} \sum_{i=0}^{n-1} \log^+
\norm {d_{p_i} f} \le \Gamma.$\fin

The proof of the theorem \ref{appli} basically
consists in producing repulsive cycles by Briend-Duval's method 
taking into account the information on inverse branches given by our
theorem \ref{nt}. \\

Let $0 < \epsi \ll \chi_1$. We introduce :
\[ \Rep_n^\epsi := \{ p \in \Rep_n \ , \ \abs {\varphi_n (p) - \Sigma_s }
\leq 2\epsi \} \] 
and write ${1 \over d_t^n} \sum_{p \in \Rep_n} \varphi_n (p) - \Sigma_s = {1 \over d_t^n} ( u_n + v_n + w_n)$, where
\[  u_n :=  \sum_{p \in \Rep_n^\epsi} \left( \varphi_n (p) -
\Sigma_s \right) \ \ , \ \   v_n :=   \sum_{p \in
  \Rep_n \setminus \Rep_n^\epsi} \left( \varphi_n (p) - \Sigma_s \right)  \
  \ , \ \   w_n :=  ( {\cd(\Rep_n) - d_t^n} )
  \Sigma_s. \] 

We show that for $n$ sufficiently large the sequences $u_n$, $v_n$ and
  $w_n$ are essentially bounded by $\epsi d_t^n$. The key estimate is given by
the following lemma whose proof is postponed to the end of the subsection.

\begin{lem}\label{l3}
There exists $n_1 \geq 1$ such that : 
\[ \forall n \geq n_1 \ , \ \cd(\Rep_n^\epsi)  \geq  d_t^n (1-\epsi)^3. \]
\end{lem}

Let us now give the expected bounds on $u_n$, $v_n$, $w_n$.
As $\Rep_n \subset \Fix_n$ and $\cd(\Fix_n) \leq (d_t^n d^n - 1)(d^n-1)^{-1}$ (recall that $d_t = d^k$), there exists $n_2 \geq 1$ such that :
\begin{equation}\label{l1} 
\forall n \geq n_2 \ , \ \cd(\Rep_n)   \leq  d_t^n (1+\epsi).
\end{equation}
Using (\ref{l1}) and the lemma \ref{l3} we see that for any $n \geq \max \{ n_1,n_2 \}$ :
\[  \abs{u_n} \leq    \cd(\Rep_n^\epsi) 2\epsi  \leq 2\epsi(1+\epsi)d_t^n \leq 4 \epsi d_t^n\ \ \textrm{ and } \ \  \abs{w_n} \leq  \ 4\Sigma_s \cdot \epsi d_t^n. \]
The lemmas \ref{l2}, \ref{l3} and the estimate (\ref{l1}) finally give :
\[   \abs{v_n} \leq   (\Gamma + \Sigma_s)  \left( \cd(\Rep_n) - \cd(\Rep_n^\epsi) \right) \leq (\Gamma + \Sigma_s)  [ (1+\epsi)-(1-\epsi)^3 ]  d_t^n, \]
which is bounded by  $5(\Gamma + \Sigma_s) \cdot \epsi d_t^n$.\\ 

We now end the proof of the theorem \ref{appli} by establishing the lemma \ref{l3}.\\

The functions $r_\epsi, T_\epsi$ and $L_\epsi$ have been introduced in the
theorem \ref{nt}. We note them shortly $r,T$ and $L$. In the sequel, 
we will not use their $\epsi$-slow/fast properties. For any $\hat x \in X$, let $n(\hat x)$ be the smallest integer satisfying
$\log L(\hat x) \leq  n \epsi$ and 
$\log T(\hat x) \leq n \epsi$. Let $n_1 \geq 1$ large enough and $r_0 > 0$ small enough such that the set
\[ \widehat H := \{ \hat x \in X \ , \  r(\hat x) \geq r_0   \ , \ n(\hat x) \leq n_1 \}\]
satisfies $\nu (\widehat H) \geq 1-\epsi/2$. We have in particular for all $\hat x
\in \widehat H$ and $n \geq n_1$ (see the theorem \ref{nt}) : 
\begin{itemize}
\item[($a_1$)] $\Lip f_{\hat x}^{-n} \leq e^{-n \chi_1 + 2n \epsi}$
  on $B_{x_0}(r_0)$,
\item[($a_2$)] for all $z \in f_{\hat x}^{-n}(B_{x_0}(r_0)), \,
  \vert  {1 \over n} \log \norm {\bigwedge^s d_z f^n} - \Sigma_s  \vert
  \leq 2\epsilon$.
\end{itemize} 

We consider two concentric families of balls $(B_i)_{1 \leq i \leq m}\subset (B_i^\ga)_{1 \leq i \leq
  m} \subset \Pj^k$ whose radii are respectively equal to $r$ and $r + \gamma \leq r_0 / 2$ (with $0 < \gamma \ll r$) and satisfy the following properties :
\begin{itemize}
\item[($b_1$)] the balls $B_i^\ga$ are disjoint,
\item[($b_2$)] $\mu (\cup_{i=1}^m B_i) \geq 1-\epsi/2$,
\item[($b_3$)] $\mu(B_i) \geq (1-\epsi) \mu (B_i^\ga)$. 
\end{itemize}
We may increase $n_1$ such that for all $n \geq n_1$ and $1 \leq i \leq m$ : 
\begin{itemize}
\item[($c_1$)] $e^{-n \lambda_1 + 2n \epsi} (r + \ga ) < \ga$ (in particular $e^{-n \lambda_1 + 2n \epsi} < 1$), 
\item[($c_2$)] if $\widehat{B_i} := \pi^{-1} (B_i)$, then $\nu \left(
  {\widehat f}^{-n} (\widehat{B_i} \cap \widehat H ) \cap \widehat{B_i} \right) \geq (1-\epsi) \nu(\widehat{B_i} \cap \widehat H) \nu(\widehat{B_i})$.
\end{itemize} 
The last assertion is a consequence of the mixing property of $\nu$.\\

Let us temporarily
fix $i \in \{1,\ldots, m\}$ and note $B:=B_i$, $B^\ga := B_i^\ga$ and
$\widehat{B^\ga} := \pi^{-1} (B^\ga)$. Observe that if $\hat x \in
\widehat H  \cap \widehat{B^\ga}$, then the set $f_{\hat
  x}^{-n}(B^\ga)$ is well defined since $\hat x \in \widehat H
\subset \{ r(\hat x)
\geq r_0 \}$ and the radius of $B^\ga$ satisfies $r + \ga \leq r_0 /2$.\\

Let us consider the collection  $\CC_n(B)$ of sets of the form $f_{\hat
  x}^{-n}(B^\ga)$ which do intersect $B$, where $\hat x \in \widehat H \cap \widehat{B^\ga}$. 
 Note that the elements of $\CC_n(B)$ are \emph{disjoint} open subsets of $\Pj^k$. 
We will establish the following estimates for $n\geq n_1$ :

\begin{equation}\label{bdcra}
\cd(\Rep_n^\epsi \cap B^\ga)  \geq  \cd(\CC_n(B))\geq d_t^n (1-\epsi)^2 \nu( \widehat H \cap
\widehat{B}).
\end{equation}

Let  $\EE$ be an element of $\CC_n(B)$ : $\EE = f_{\hat
  x}^{-n}(B^\ga)$ intersects $B$ and $\hat x \in \widehat H \cap
\widehat{B^\ga}$. By ($a_1$) and ($c_1$), $f_{\hat
  x}^{-n}$ contracts on $B^\ga \subset B_{x_0}(r_0)$. Moreover $\EE$ is
contained in $B^\ga$, because $\EE$  intersects $B$ and its diameter is
less than $\ga$ by ($c_1$). We thus obtain a point $p \in \EE$
which is fixed by $f_{\hat x}^{-n} : B^\ga \to B^\ga$. This point is
$n$-periodic and repulsive for $f$, so $p \in
\Rep_n$. We have also $p \in  f_{\hat x}^{-n}(B^\ga) \subset  f_{\hat
  x}^{-n}(B_{x_0}(r_0))$ which leads to $p \in
\Rep_n^\epsi$ by ($a_2$). This gives the first inequality in (\ref{bdcra}).

 Let us now justify the second inequality.
According to ($c_2$) and the relation $\mu = \pi_* \nu$ we obtain :
\[  (1-\epsi) \nu(\widehat{B} \cap \widehat H) \mu(B) \leq \mu \left[
  \pi \left( {\widehat f}^{-n} (\widehat H  \cap \widehat B)  \right)
  \cap B  \right] \leq \mu  \Big( \bigcup f_{\hat x}^{-n}(B^\ga) \Big), \]
where the union runs over the elements of $\CC_n(B)$. By the jacobian property we have $\mu (f_{\hat
  x}^{-n}(B^\ga)) = \mu(B^\ga) / d_t^n$ and the right hand side is thus equal to $\cd(\CC_n(B)) \mu(B^\ga) / d_t^n$. We then get the desired inequality by ($b_3$).\\

The estimates (\ref{bdcra}) imply finally :
\begin{equation}\label{fgh}
  \cd(\Rep_n^\epsi) \geq \sum_{i=1}^m  \cd(\Rep_n^\epsi \cap B^\ga_i)
\geq \sum_{i=1}^m   \cd(\CC_n(B_i)) \geq d_t^n (1-\epsi)^2 \nu( \widehat H \cap
\cup_{i=1}^m \widehat{B_i}). 
\end{equation}
 As $\nu
(\widehat H) \geq 1-\epsi/2$ and $\nu( \cup_{i=1}^m \widehat{B_i}) = \mu ( \cup_{i=1}^m B_i ) \geq 1 - \epsi / 2$ (see ($b_2$)) we have $\nu( \widehat H \cap
\cup_{i=1}^m \widehat{B_i}) \geq 1-\epsi$. So, by (\ref{fgh}),
$\cd(\Rep_n^\epsi) \geq d_t^n (1-\epsi)^3$, and the lemma
\ref{l3} is proved.

\section{Appendix : properties of resonant maps}\label{pest}

We prove here the proposition \ref{olk} and \ref{estreso} we stated in the subsection \ref{bresonn}.
A result similar to the proposition \ref{olk} may be found in the article \cite{GK}, lemma 1.1.

\subsection{Proof of the proposition \ref{olk}} 

We want to prove that for all $2 \leq j \leq k$, the component 
\begin{equation*}
\pi_j(K_{\tau(x)}) (K_x) = \pi_j(K_{\tau(x)}) (\pi_1 (K_x),\ldots,\pi_{j-1} (K_x))
\end{equation*}
is made of $j$-resonant monomials (the dependance of the $j-1$ first
components in the right hand side is a consequence of (\ref{lok}), subsection \ref{bresonn}). For simplicity, we may assume that $k_j =
1$ for all $1 \leq j \leq k$, therefore $\lambda_j = \Lambda_j$.  As
we are concerned with a degree property, we do not compute the
coefficients of polynomials (we set them equal to $1$). Recall that $P_{\alpha}=z_1^{\alpha_1}\cdot\cdot\cdot z_k^{\alpha_k}$.
By writing $\pi_j(K_{\tau(x)}) = \sum_{\beta \in \RRR_j} P_\beta$, 
it suffices to consider the following polynomial, where $\beta \in \RRR_j$ : 
\begin{equation}\label{expa1}
P_\beta (\pi_1 (K_x),\ldots,\pi_{j-1} (K_x)) = (\pi_1 (K_x))^{\beta_1} \ldots (\pi_{j-1} (K_x))^{\beta_{j-1}}. 
\end{equation}
For all $1 \leq l \leq j-1$, we write $\pi_l (K_x) =
\sum_{i=0}^{\cd(\RRR_l)} P_{\al(l,i)}$, where $\al(l,0) =
(0,\ldots,1,\ldots,0)$ (so that $P_{\al(l,0)}=z_l$) and $\{ \al(l,i)
\, , \, 1 \leq i \leq  \cd(\RRR_l) \} = \RRR_l$. Observe that we have for any $1 \leq l \leq j-1$ and $0 \leq i \leq  \cd(\RRR_l)$ :
\begin{equation}\label{alph}
\sum_{p=1}^{j-1} \al_p(l,i) \Lambda_p = \sum_{p=1}^l \al_p(l,i) \Lambda_p  =  \Lambda_l.
\end{equation}
We now expand (\ref{expa1}) and prove that we get a sum of $j$-resonant monomials. By the Newton's formula, the monomials 
that appear in the expansion of (\ref{expa1}) have the form : 
\begin{equation}\label{expa2}
P_\ga = \prod_{i=0}^{\cd(\RRR_1)} P_{\al(1,i)}^{\beta_1(i)} \ \ \ldots  \prod_{i=0}^{\cd(\RRR_{j-1})} P_{\al(j-1,i)}^{\beta_{j-1}(i)}  \ ,  
\end{equation}
where $\sum_{i=0}^{\cd(\RRR_l)} \beta_l(i) = \beta_l$ for all $1 \leq
l \leq j-1$. We have to prove that $\ga \in \RRR_j$, that is : $\sum_{p=1}^{j-1}  \ga_p \Lambda_p = \Lambda_j$. First observe that $\abs{\ga} \geq \abs{\beta} \geq 2$ and that $P_\ga = z_1^{\ga_1}\ldots
z_{j-1}^{\ga_{j-1}}$ (indeed there are no $z_j, \ldots , z_k$ in the right hand side of (\ref{expa1})). By  (\ref{expa2}) we have for all $1 \leq p \leq j-1$ :
\begin{equation*}
\ga_p = \sum_{i=0}^{\cd(\RRR_1)} \al_p(1,i)\beta_1(i) +  \ldots + \sum_{i=0}^{\cd(\RRR_{j-1})} \al_p(j-1,i)\beta_{j-1}(i) = 
\sum_{l=1}^{j-1} \sum_{i=0}^{\cd(\RRR_l)} \al_p(l,i)\beta_l(i). 
\end{equation*}
Thus the sum $\sum_{p=1}^{j-1}  \ga_p \Lambda_p$ is equal to :
\begin{equation}\label{expa3}
     \sum_{p=1}^{j-1} \Lambda_p \left( \sum_{l=1}^{j-1} \sum_{i=0}^{\cd(\RRR_l)} \al_p(l,i)\beta_l(i) \right)  
= \sum_{l=1}^{j-1}  \left(  \sum_{i=0}^{\cd(\RRR_l)} \beta_l(i) \left( \sum_{p=1}^{j-1} \al_p(l,i) \Lambda_p \right) \right).                                    
\end{equation}
By using successively (\ref{alph}), $\sum_{i=0}^{\cd(\RRR_l)} \beta_l(i) = \beta_l$ and $\beta \in \RRR_j$, the right hand of (\ref{expa3}) is equal to 
\[ \sum_{l=1}^{j-1}  \left( \sum_{i=0}^{\cd(\RRR_l)} \beta_l(i) \, \Lambda_l \right)  =  \sum_{l=1}^{j-1}  \beta_l \Lambda_l = \Lambda_j. \]
We have therefore $\ga \in \RRR_j$, which completes the proof of the proposition.

\subsection{Proof of the proposition \ref{estreso}} 

Let $\RR$ be a resonant bundle map. By the proposition \ref{olk},
the degree of the iterates of $\RR$ is bounded by
${\tilde q}$. We define the function $\norm{\RR} :=  \max \{  \abs
{\RR^{(1)}} , \ldots , \abs{\RR^{({\tilde q})}} \}$ : for all $x \in X$,
$\norm{\RR}(x)$ is the maximum of the coefficients of the polynomial
map $R_x$ (see subsection
\ref{tube} for the definition of $\abs{\RR^{(j)}}$). We show in the next proposition that for all $1 \leq j \leq l$, $\norm{\pi_j(\RR^n)}$ is close to $\norm{\pi_j(\AA^n)}$. This fact will be useful for proving the proposition \ref{estreso}.

\begin{prop}\label{pj}
Let $\RR : E(\rho_\epsi) \to E(\rho_\epsi)$ be a resonant bundle map where $\rho_\epsi$ is a $\epsi$-slow function. There exist
$\theta \geq 1$ (depending only on $\tilde q$) and a $\theta\epsi$-fast
function $M_{\theta\epsi} : X \to [1,+\infty[$ such that 
\[  \forall j \in \{ 1,\ldots,l \} \ , \  \forall n \geq 1 \ , \ e^{n(\Lambda_j - \epsilon)} \leq  \norm{ \pi_j(\RR^n) } \leq M_{\theta\epsi} e^{n(\Lambda_j + \theta \epsilon)}.  \] 
\end{prop}

\proof The estimate from below follows from  $\norm{ \pi_j(\RR^n) }
\geq \abs { \pi_j(\AA^n) } \geq e^{n(\Lambda_j - \epsilon)}$. We prove
the estimate from above. We may assume that $k_j = 1$ and
$\lambda_j = \Lambda_j$ for all $1 \leq j \leq k$. Let also $\Delta :=
\cd(\RRR) + 1$. Let $1 \leq s \leq k$ and define the assertion $(i_s)$ : 
there exists a $(\tilde q + \ldots + {\tilde q}^{s-1})\epsi$-fast function $M_s : X \to [1,+\infty[$ such that 
\begin{equation*}
 \forall j \in \{ 1,\ldots,s \} \ , \  \forall n \geq 1  \ , \  \norm{ \pi_j(\RR^n) } \leq M_s e^{n(\Lambda_j + {\tilde q}^{2s} \epsilon)} .  
\end{equation*}
Let  $\psi_{\tilde q \epsi}$ be a $\tilde q \epsi$-fast function such that $\psi_{\tilde q \epsi} \geq \norm{\RR}$ (cf lemma
\ref{cow}). We proceed by induction on $s$. The assertion $(i_1)$ is satisfied with $M_1 := 1$. Indeed,
$\pi_1(\RR^n) = \pi_1(\AA^n)$ because $\RRR_1$ is empty (see subsection \ref{bresonn}). 
Assume that $(i_{s-1})$ is satisfied for $1 \leq
s-1 \leq k-1$. We define
\[  M_{s} :=  {\psi_{\tilde q \epsi}  M_{s-1}^{{\tilde q}} \Delta^{{\tilde q}+1} e^{-\Lambda_{s}}  \over 1 -
  e^{-( {\tilde q}^{2s} -1)\epsi} } . \] 
Observe that $M_s$ is $(\tilde q + \ldots +{\tilde q}^{s-1})\epsi$-fast. As $M_{s} \geq M_{s-1}$, the assertion $(i_s)$ is proved if we
establish the following assertions for all $n \geq 1$ :   
\[ (ii_n) : \norm{ \pi_{s} (\RR^n) } \leq M_s e^{n(\Lambda_{s} + {\tilde q}^{2s}\epsilon)}.  \]
The assertion $(ii_1)$ is fulfilled because :
\[ \norm { \pi_{s}(\RR)} \leq \norm{\RR} \leq \psi_{\tilde q \epsi} \leq \psi_{\tilde q \epsi} {M_{s-1}^{{\tilde q}} \Delta^{{\tilde q}+1} e^{{\tilde q}^{2s}\epsi}  \over 1 -
  e^{-( {\tilde q}^{2s} -1)\epsi} }  = M_s e^{\Lambda_{s} + {\tilde q}^{2s}\epsilon}.  \]
Let us assume that $(ii_n)$ is
  true and establish $(ii_{n+1})$. We have  
\begin{equation}\label{bbn}
 \pi_{s} (R_x^{n+1}) =  \pi_{s} (R_{\tau^n(x)}) ( \pi_1 (R_x^n),\ldots,\pi_{s} (R_x^n)). 
\end{equation}
By (\ref{lok}) (subsection \ref{bresonn}), the polynomial $\pi_{s} (R_{\tau^n(x)})$ has the form \[ \pi_{s} (R_{\tau^n(x)}) (z_1,\ldots,z_{s}) = \si(\tau^n(x)) z_{s} + \sum_{\al \in \RRR_{s}} \si_{\al}(\tau^n(x)) P_\al(z_1,\ldots,z_{s-1}), \]
where $\si(\tau^n(x)),\si_{\al}(\tau^n(x)) \in \C$. So by (\ref{bbn}), $\pi_{s} (R_x^{n+1})$ is the sum of $I_{n}$ and $J_{n}$ : 
\[ I_{n} := \si(\tau^n(x)) \pi_{s}(R_x^n), \] 
\begin{equation}\label{ttt}
 J_{n} := \sum_{\al \in \RRR_{s}} \si_{\al}(\tau^n(x)) P_\al(\pi_1 (R_x^n),\ldots,\pi_{s-1}(R_x^n) ).
\end{equation}
To simplify the exposition, we note in the sequel $u := e^{\Lambda_{s} + {\tilde q}^{2s}
  \epsilon}$ and $M_s$ for $M_s(x)$, $\psi_{\tilde q \epsi}$ for $\psi_{\tilde q \epsi}(x)$, etc.  We obtain by  $(ii_n)$ and $\abs{\si(\tau^n(x))} \leq e^{\Lambda_{s} + \epsi}$ :
\[ \norm{I_{n}}  \leq  e^{\Lambda_{s} + \epsi} M_{s}  u^n = \left[
M_{s} e^{-({\tilde q}^{2s} -1)\epsi} \right]  u^{n+1}. \]
The end of the proof consists in verifying the following estimate, by using $(i_{s-1})$ : 
\begin{equation}\label{efg}
\norm{J_{n}}  \leq \left[  \Delta^{{\tilde q}+1} \psi_{\tilde q \epsi}   M_{s-1}^{{\tilde q}} e^{-\Lambda_{s}} \right] u^{n+1}.\end{equation}
Indeed the two preceding lines imply $(ii_{n+1})$ by the very definition of $M_{s}$ :
\[ \norm{ \pi_{s} (\RR^{n+1}) }
  \leq \norm{I_{n}} + \norm{J_{n}} \leq \left[   M_{s} e^{-({\tilde q}^{2s} -1)\epsi} +  \Delta^{{\tilde q}+1} \psi_{\tilde q \epsi}  M_{s-1}^{{\tilde q}} e^{-\Lambda_{s}}  \right]  u^{n+1} =  M_{s} u^{n+1}. \]
We now prove (\ref{efg}). By the proposition \ref{olk}, $J_n$ is a linear combination of resonant monomials.
Recall also that $\norm{J_{n}}$ is the maximum of
the coefficients of $J_{n}$ in the basis $(P_\beta)_{\beta \in
  \RRR_{s}}$. Fix a resonant degree $\beta \in
\RRR_{s}$. Let $\al \in \RRR_s$ and expand $P_\al(\pi_1
(R_x^n),\ldots,\pi_{s-1}(R_x^n))$. Observe that for all
$1 \leq j \leq s-1$, $\pi_j(R^n_x)$ is a sum of at most $\Delta = \cd(\RRR)+1$
monomials, because $\RR$ is resonant. So the degree $\beta\in \RRR_s$ appears in the preceding expansion at most 
$\Delta^{\al_1 + \ldots +\al_{s-1}}$ times. It implies that the $P_\beta$-coefficient $\si_{\beta}(\al)$ of 
$\si_{\al}(\tau^n(x)) P_\al(\pi_1 (R_x^n),\ldots,\pi_{s-1}(R_x^n))$ satisfies by $(i_{s-1})$ and $\abs{\si_{\al}(\tau^n(x))} \leq \psi_{\tilde q \epsi}(x) e^{n\epsi}$ :
\begin{displaymath}
\begin{array}{rcl}
   \abs{\si_{\beta}(\al)}  &  \leq   & \abs{\si_{\al}(\tau^n(x))} \Delta^{\al_1 + 
\ldots +\al_{s-1}} \norm{\pi_1(R_x^n)}^{\al_1} \ldots \norm{\pi_{s-1}(R_x^n)}^{\al_{s-1}}     \\
       &  \leq   &   \psi_{\tilde q \epsi}(x) e^{n\epsi}  
[ \Delta M_{s-1}  
e^{n {\tilde q}^{2(s-1)} \epsilon}] ^ { \al_1 + \cdots + \al_{s-1}}  e^{n ( \al_1 \Lambda_1 + \cdots + \al_{s-1} \Lambda_{s-1}  )}  .
\end{array}                                      
\end{displaymath}
Observe that $\al \in \RRR_{s}$ implies $\al_1 + \ldots +\al_{s-1} \leq {\tilde q}$ and 
$\al_1 \Lambda_1 + \ldots + \al_{s-1} \Lambda_{s-1} = \Lambda_{s}$. We deduce 
(use $1 + {\tilde q}{\tilde q}^{2(s-1)} \leq {\tilde q}^2 {\tilde q}^{2(s-1)} = {\tilde q}^{2s}$) :
\[ \abs{ {\si_{\beta}}(\al) } \leq  [ \psi_{\tilde q \epsi}  \Delta^{{\tilde q}}   M_{s-1}^{\tilde q} e^{-\Lambda_s} ] 
e^{(n+1) {\tilde q}^{2s} \epsilon}  e^{(n+1) \Lambda_{s} } = [ \psi_{\tilde q \epsi} \Delta^{{\tilde q}}   M_{s-1}^{\tilde q} e^{-\Lambda_s} ] u^{n+1}.  \]
We now use (\ref{ttt}) to obtain $\norm{J_n} \leq \sum_{\al \in \RRR_s} \abs{
  {\si_{\beta}}(\al) } \leq \Delta [\psi_{\tilde q \epsi} \Delta^{{\tilde q}}
  M_{s-1}^{\tilde q} e^{-\Lambda_s} ] u^{n+1}$, this completes the proof of
  (\ref{efg}). We let finally $\theta = \max \{ \tilde q + \ldots + {\tilde q}^{k-1} ,  {\tilde
    q}^{2k} \}$ to get the two assertions : the function $M_{\theta\epsi} := M_k$ is $\theta\epsi$-fast and  
$\norm{ \pi_j(\RR^n) } \leq M_{\theta\epsi} e^{n(\Lambda_j + \theta \epsilon)}$. \fin

We now prove the proposition \ref{estreso}. Recall that 
\[ (\lambda_1, \ldots , \lambda_k) =
(\Lambda_1 , \ldots , \Lambda_1, \cdots ,  \Lambda_j , \ldots ,
\Lambda_j ,\cdots, \Lambda_l, \ldots, \Lambda_l), \]
where $\Lambda_j$ appears $k_j$ times. We note $(e_i)_i$ the canonical orthonormal basis of 
$\C^k$. If $A:\C^k\to\C^k$ is linear, we define $\bigwedge^s A:\bigwedge^s \C^k\to\bigwedge^s\C^k$ as the linear extension of the map
$L$ satisfying $L(e_{i_1}\wedge\cdot\cdot\cdot \wedge e_{i_s}):=A(e_{i_1})\wedge\cdot\cdot\cdot A(e_{i_s})$ for any $1 \le i_1<\cdot\cdot\cdot <i_s\le k$.

\begin{proposition}{\ref{estreso}}
Let $\RR : E(\rho_\epsi) \to E(\rho_\epsi)$ be a resonant bundle map
    where $\rho_\epsi$ is a $\epsi$-slow function. Let
    $\eta := k \theta$ ($\theta$ depends on $\tilde q$, see proposition \ref{pj}). There exists a $\eta\epsi$-fast function $H_{\eta\epsi} : X \to [1,+\infty[$ such that
    for all $w \in E_x(\rho_\epsi)$, $s \in \{1,\ldots,k \}$ and $n
    \geq 1$ :  $\Big \vert   {1 \over n} \log \abs{\bigwedge^s d_w
    R^n_x} - (\lambda_1 + \ldots + \lambda_s)  \Big \vert  \leq {1
    \over n} \log H_{\eta\epsi}(x) + \eta \epsilon$.
\end{proposition}

\proof We first give the proof in the case $k_j = 1$ for all $1 \leq j \leq k$ (it implies $\lambda_j = \Lambda_j$). Let $x \in X$ and $w \in E_x(\rho_\epsi)$. By the remark \ref{trian} (subsection \ref{bresonn}), the matrix of $d_w R^n_x$ in the canonical basis $(e_i)_{1 \leq i \leq k}$ is lower triangular, so we have :
\[  d_w R^n_x (e_i) = \zeta_n(i) e_i + \sum_{j=i+1}^k \om_n(i,j) e_j \ \textrm{ and } \  d_w R^n_x (e_k) = \zeta_n(k) e_k. \]
We fix $1 \leq i \leq k$ and $i +1 \leq j \leq k$. We give bounds for $\zeta_n(i)$ and $\om_n(i,j)$. We have 
\begin{equation}\label{dcb1}
e^{n(\Lambda_i - \epsi)} \leq
\abs{\zeta_n(i)} \leq  e^{n(\Lambda_i + \epsi)}
\end{equation}
because the diagonal part of $d_w R^n_x$ is equal to $A^n_x$ (see
remark \ref{trian}). Recall that $\pi_j(R^n_x)$ is a sum of at most
$\cd(\RRR_j)$ resonant monomials with degree $\alpha$ lower than ${\tilde q}$
(see the proposition \ref{olk}). We deduce by the proposition \ref{pj}
and $\rho_\epsi \leq 1$ that : 
\begin{equation}\label{dcb2}
\abs{\om_n(i,j)} \leq  \sum_{\al \in \RRR_j} \norm{\pi_j (R^n_x)} \abs{\al} \rho_\epsi(x)^{\abs{\al}-1}  \leq  \cd(\RRR_j)   
M_{\theta\epsi}(x) e^{n(\Lambda_j + \theta \epsi)} {\tilde q}.
\end{equation}
Let $M_{\theta\epsi}'(x) := \cd(\RRR) {\tilde q}  M_{\theta\epsi}(x)\ge 1$. In the sequel, we will multiply $M_{\theta\epsi}'$ 
by constants depending only on $k$ and $s$ without mentionning it. By (\ref{dcb1}), (\ref{dcb2}) and the inequality $\Lambda_j < \Lambda_i$, we get :
\begin{equation}\label{dcb3}
 \abs{d_w R^n_x (e_i)} \leq M_{\theta\epsi}'(x) e^{n(\Lambda_i + \theta \epsi)}.
\end{equation}
Now we focus on the coefficients of the matrix $\bigwedge^s d_w R^n_x$ in the orthonormal basis $e_{i_1} \wedge \ldots \wedge e_{i_s}$, 
with $1 \leq i_1 < \ldots < i_s \leq k$. 
We begin with the vector $\bigwedge^s d_w R^n_x ( e_1 \wedge \ldots \wedge e_s )$ which is equal to :
\begin{equation}\label{eee}
 \zeta_n(1) \ldots \zeta_n(s) e_1 \wedge \ldots \wedge e_s + \sum_{(i_1,\ldots,i_s) \neq (1,\ldots,s)} \om_n'(i_1,\ldots,i_s) e_{i_1} \wedge \ldots \wedge e_{i_s}.
\end{equation}
Observe now the following estimates for $(i_1,\ldots,i_s) \neq
(1,\ldots,s)$ and $n \geq 1$. The first inequality is a consequence of
(\ref{dcb2}) :
\begin{equation} \label{sxc}
\abs{ \om_n'(i_1,\ldots,i_s) } \leq
M_{\theta\epsi}'(x)^s e^{n(\Lambda_{i_1} + \ldots + \Lambda_{i_s} + s \theta
  \epsi)} \leq M_{\theta\epsi}'(x)^s e^{n(\Lambda_1 + \ldots + \Lambda_s +  k\theta\epsi)}.
\end{equation}
We obtain with the lines (\ref{dcb1}), (\ref{eee}) and (\ref{sxc}) that for any $n \geq 1$ : 
\begin{equation}
  e^{n(\Lambda_{1} + \ldots + \Lambda_{s} -s \epsi)} \leq \abs{ \bigwedge^s d_w R^n_x ( e_1 \wedge \ldots \wedge e_s )  }    \leq M_{\theta\epsi}'(x)^k e^{n(\Lambda_{1} + \ldots + \Lambda_{s} + k\theta\epsi)}.
\end{equation}
We now focus on $\bigwedge^s d_w R^n_x (e_{i_1} \wedge \ldots \wedge e_{i_s})$. 
The Hadamard's inequality implies with (\ref{dcb3}) :
\begin{displaymath}\label{pp}
\begin{array}{rcl}
    \abs{ \bigwedge^s d_w R^n_x (e_{i_1} \wedge \ldots \wedge e_{i_s}) }    &  \leq  & \abs{ d_w R^n_x (e_{i_1}) } \ldots \abs{ d_w R^n_x (e_{i_s})}  \\
       & \leq  &  M_{\theta\epsi}'(x)^s e^{n(\Lambda_{i_1} + \ldots + \Lambda_{i_s} + k \theta \epsi)} \\
       & \leq  &  M_{\theta\epsi}'(x)^k  e^{n(\Lambda_{1} + \ldots +
       \Lambda_{s} + k \theta \epsi)}. 
\end{array}                                      
\end{displaymath} 
We finally obtain the following estimates for any $n \geq 1$, where $\eta = k\theta$ and $H_{\eta\epsi}$ is equal to the $ \eta\epsi$-fast function $(M_{\theta\epsi}')^k$ :
\[ \Big \vert   {1 \over n} \log \abs{\bigwedge^s d_w R^n_x} -
(\lambda_1 + \ldots + \lambda_s)  \Big \vert  \leq {1 \over n} \log
H_{\eta\epsi}(x) + \eta \epsilon.\]

We sketch the proof in the general case, i.e. when $k = \sum_{i=1}^l
k_j$ and $k_j \geq 1$. We may assume that the block diagonal matrix $A_x^n$ is a 
lower (block) triangular matrix. Indeed, for all $1 \leq j \leq l$, there exists a
matrix $V_j  \in \U_{k_j}(\C)$ such that $T_{j}^n := V_j [
  \pi_j(A^n_x) ] V_j^{-1}$ is lower triangular. The metric property :
\[ \forall (x,v) \in \LL^j \ , \ e^{n(\Lambda_j - \epsilon)} \abs{v}
\leq \abs{\pi_j (A^n_x)(v)}  \leq e^{n(\Lambda_j + \epsilon)}  \abs{v} \]
implies that the modulus of the coefficients of $T_{j}^n$ are $\leq e^{n(\Lambda_j + \epsilon)}$, and that the diagonal coefficients $(\zeta_n^{(j)}(i))_{1\leq i \leq k}$ of  $T_{j}^n$  satisfy :
\begin{equation*}
e^{n(\Lambda_j - \epsi)} \leq
\abs{\zeta_n^{(j)}(i)} \leq  e^{n(\Lambda_j + \epsi)}.
\end{equation*}     
We thus obtain estimates analogous to (\ref{dcb1}). Let $V$ be the block
diagonal matrix $(V_1,\ldots,V_l)$. The matrix $T_x^n := V
[A^n_x] V^{-1}$ is therefore block diagonal, with lower triangular blocks
$(T_{1}^n,\ldots,T_{l}^n)$. The matrix $S_x^n :=
V [ d_w R^n_x ] V^{-1}$ is also lower triangular (see remark \ref{trian}) and its coefficients outside the block diagonal matrix $T_x^n$ satisfy estimates analogous to (\ref{dcb2}) (the $(V_j)_{1\leq j \leq l}$ are unitary transformations). We deduce as before that for any $n \geq 1$ :
\[ \Big \vert   {1 \over n} \log \abs{\bigwedge^s S_x^n} - (\lambda_1
+ \ldots + \lambda_s)  \Big \vert  \leq {1 \over n} \log
H_{\eta\epsi}(x) + \eta \epsilon. \]
The conclusion follows from  $\abs{\bigwedge^s S_x^n} =
\abs{\bigwedge^s V \, \left[ \bigwedge^s d_w R^n_x \right] \,
  \bigwedge^s V^{-1}  } = \abs{\bigwedge^s d_w R^n_x}$. \fin

\vspace{1 cm}

{\footnotesize F. Berteloot}\\
{\footnotesize Universit\'e Toulouse III}\\
{\footnotesize Institut Mathématique de Toulouse}\\
{\footnotesize Equipe Emile Picard, Bat. 1R2}\\
{\footnotesize 118, route de Narbonne }\\
{\footnotesize F-31062 Toulouse Cedex 9, France}\\
{\footnotesize berteloo@picard.ups-tlse.fr}\\ 

{\footnotesize C. Dupont}\\
{\footnotesize Universit\'e Paris XI-Orsay}\\
{\footnotesize CNRS UMR 8628}\\
{\footnotesize Math\'ematique, B\^at. 425}\\
{\footnotesize F-91405 Orsay Cedex, France}\\
{\footnotesize christophe.dupont@math.u-psud.fr}\\

{\footnotesize L. Molino}\\
{\footnotesize Universit\`a di Parma}\\
{\footnotesize Dipartimento di Matematica}\\
{\footnotesize Parco Area delle Scienze, Viale Usberti 53/A}\\
{\footnotesize I-43100 Parma, Italia.}\\
{\footnotesize laura.molino@unipr.it}

\end{document}